\documentclass[preprint]{elsarticle}

\journal{Journal of Computational Physics}

\usepackage{tikz}
\usepackage[inkscapelatex=false]{svg}
\usepackage{mathdots}
\usepackage{yhmath}
\usepackage{cancel}
\usepackage{color}
\usepackage{siunitx}
\usepackage{array}
\usepackage{multirow}
\usepackage{multicol}
\usepackage{amssymb}
\usepackage{tabularx}
\usepackage{subfigure}
\usepackage{extarrows}
\usepackage{booktabs}
\usetikzlibrary{fadings}
\usetikzlibrary{patterns}
\usetikzlibrary{shadows.blur}
\usetikzlibrary{shapes}
\usepackage{amsmath}
\usepackage[utf8]{inputenc} 
\usepackage[T1]{fontenc}    
\usepackage{hyperref}       
\usepackage{url}            
\usepackage{booktabs}       
\usepackage{amsfonts}       
\usepackage{nicefrac}       
\usepackage{microtype}      
\usepackage{lipsum}		
\usepackage{graphicx}
\usepackage{doi}
\usepackage{comment}
\usepackage{pgf}

\newcommand{\V}[1]{\textbf{#1}}                 

\newcommand{\Partial}[2]{\frac{\partial #1}{\partial #2}}

\newcommand{\myblue}[1]{{\color[rgb]{0.0,0.0,0.65} #1}}
\newcommand{\myred}[1]{{\color[rgb]{0.65,0.0,0.0} #1}}

\begin{document}

\begin{frontmatter}

\title{A particle-in-Fourier method with \myblue{semi-discrete} energy conservation for non-periodic boundary conditions}
 
 \author[Nigel]{Changxiao Nigel Shen}
 \author[Antoine]{Antoine Cerfon}
 \author[sriram]{Sriramkrishnan Muralikrishnan}
 \address[Nigel]{Massachussets Institute of Technology, USA}
 \address[Antoine]{Type One Energy Group - Canada Inc., Canada}
 \address[sriram]{Jülich Supercomputing Centre, Forschungszentrum Jülich GmbH, Germany}

\begin{abstract}

We introduce a novel particle-in-Fourier (PIF) scheme based on \cite{Evstatiev_2013,mitchell2019efficient} that extends its applicability to non-periodic boundary conditions. Our method handles free space boundary conditions by replacing the Fourier Laplacian operator in PIF with a mollified Green's function as first introduced by Vico-Greengard-Ferrando \cite{Vico_2016}. This modification yields highly accurate free space solutions to the Vlasov-Poisson system, while still maintaining energy conservation up to an error bounded by the time step size. We also explain how to extend our scheme to arbitrary Dirichlet boundary conditions via standard potential theory, which we illustrate in detail for Dirichlet boundary conditions on a circular boundary. We support our approach with proof-of-concept numerical results from two-dimensional plasma test cases to demonstrate the accuracy, efficiency, and conservation properties of the scheme. \myred{By avoiding grid heating and finite grid instability we are able to show an order of magnitude speedup compared to the standard PIC scheme for a long time integration cyclotron simulation.}

\end{abstract}


\end{frontmatter}

\section{Introduction}
The particle-in-cell (PIC) method has become a common approach to numerically solve kinetic equations in plasma physics \cite{birdsall1969clouds, buneman1963computer, dawson1962one, langdon1970theory, morse1970multidimensional}. It is physically intuitive, comparatively easy to parallelize, and robust for a wide range of physical scenarios. Traditional PIC schemes, however, are known to be unable to preserve energy \cite{birdsall2018plasma, langdon1970effects}, leading to a non-physical phenomenon known as "grid-heating" and numerical instability for long-time simulations \cite{okuda1972nonphysical}. Earlier attempts to improve the conservation properties of PIC fall into several categories: some are based on the Lagrangian formulation in \cite{low1958lagrangian}, such as \cite{eastwood1991virtual, lewis1970energy}; others take advantages of implicit integrators, such as \cite{chen2011energy, markidis2011energy}. Structure preserving geometric PIC schemes are also being developed based on a variational formulation \cite{campos2022variational,squire2012geometric} and discretization of the underlying Hamiltonian structure \cite{kraus2017gempic,he2016hamiltonian}. Although many of these schemes have excellent long time stability and conservation properties, they fail to conserve one of either charge or momentum. They also differ significantly from the standard PIC framework and hence may not be as intuitive or easy to transition to from an application and implementation point of view.

Among recent developments for structure-preserving PIC, the particle-in-Fourier (PIF) scheme developed in \cite{Evstatiev_2013, mitchell2019efficient} retains the intuitive nature of standard PIC schemes as well as their relative simplicity, and at the same time has excellent conservation properties. By exploiting the non-uniform fast Fourier transform \cite{DuttRokhlin, Beylkin1995OnTF, press1989fast}, PIF directly evaluates the Fourier modes of the field efficiently, without any spatial charge spreading and interpolation processes. Both charge and momentum are conserved in this scheme, and energy is conserved in the continuous time limit. The finite-grid instability and mode-coupling effects are also eliminated \cite{huang2016finite}. \myred{In \cite{Evstatiev_2013} the authors show that in finite dimensions PIF is the only scheme which can conserve both energy and momentum simultaneously before time discretization. Even though the exact energy conservation is lost after time discretization with explicit time integrators, the scheme still possesses excellent stability and convergence properties for long time integration simulations as shown in Section \ref{grid_heating}. PIF scheme combined with implicit time integration has not been explored yet and could be an interesting research direction to explore in the future.} The use of the non-uniform fast Fourier transform reduces the complexity of the algorithm from $\mathcal{O}(N_pN_m^d)$ to $\mathcal{O}(N_p+N_m^d \log N_m^d)$, where $N_p$ is the total number of particles and $N_m$ is the number of Fourier modes used in the simulation. Such a Fourier approach has also proved to be successful in other numerical contexts, such as in its application for the immersed boundary method \cite{chen2023fourier}.

On the other hand, the choice of a Fourier basis is naturally incompatible with non-periodic boundary conditions, such as free space and Dirichlet boundary conditions. Although particle methods are frequently applied to periodic problems \cite{ji2022magnetic, Klion_2023}, most problems of physics and engineering interest cannot be approximated accurately with periodic boundary conditions \cite{ji2022magnetic,LeBars2022,Taccogna2023}. One of the strengths of particle methods have been their adaptability to complex geometries and boundary conditions \cite{Ku2018,Cole2019,LeBars2022}. Limiting PIF to periodic problems would severely limit its range of applicability and its desirability. To handle nonperiodic boundary conditions, one could consider other basis functions, such as Legendre polynomials. Unfortunately, no basis other than the Fourier basis has been found that preserves the spatial translational invariance of the Lagrangian describing the kinetic plasma of interest  \cite{Evstatiev_2013,ameres2018stochastic}.

In this paper, we propose a modified version of PIF scheme enabling the method to adequately handle free space boundary conditions and Dirichlet boundary conditions. We replace the periodic Fourier Laplacian operator $1/k^2$ used in PIF with a spatially truncated free space Green's function $\hat{g}^L(\textbf{k})$ proposed in \cite{Vico_2016}, with its truncation radius $L$ determined based on our computational domain and the shape function that we choose. It has been shown in \cite{Vico_2016, zou2021fftbased} that the Fourier kernel $\hat{g}^L$ gives us a spectrally accurate solution to the free space Poisson problem for any smooth source $\rho(\textbf{x})$ on a Cartesian grid. We show that this kernel yields a spectrally accurate solution even if the source comes from Lagrangian particles, which are not aligned with grid points. Therefore, by applying $\hat{g}^L$ to our PIF scheme, we are able to obtain free space solutions to the Vlasov-Poisson system with high accuracy and exact charge, momentum and semi-discrete energy conservation. Additionally, by combining the free space PIF algorithm with a Laplace solver based on boundary integral method, we can obtain the solution to any Dirichlet problem.

The remainder of this paper is structured as follows. In Section \ref{sec:pre}, we revisit the numerical details of the particle-in-cell and particle-in-Fourier algorithms; then we give a brief overview of the free space Poisson solver proposed by Vico, Greengard and Ferrando \cite{Vico_2016}, which plays a key role in our new numerical scheme. In Section \ref{sec:combine}, we present our strategy of combining the free space Poisson solver and the Laplace solver with the PIF algorithm, followed by an error analysis for the conservation of energy and momentum. In Section \ref{sec:tests}, we report our numerical results using two-dimensional test cases. In Section \ref{sec:con}, we summarize our work and suggest directions for future work.
\section{Prerequisites}
\label{sec:pre}
Among its various applications, PIC is very often used to solve the kinetic equations that arise in plasma physics \cite{hockney1988computer, birdsall2018plasma}. In this work, we will focus on the kinetic equation describing the evolution of a single species, collisionless, electrostatic plasma immersed in a constant, uniform magnetic field $\textbf{B}_0$. The evolution of the distribution function $f(\textbf{x}, \textbf{v}, t)$ of that single-species plasma is described by the following modified Vlasov-Poisson system (here we set the vacuum permittivity $\epsilon_0$ to 1 for simplicity):
\begin{align}\Partial{f}{t}+\V{v}\cdot\nabla_{\V{x}} f+\frac{q}{m}(\V{E}+\V{v}\times \textbf{B}_0)\cdot\nabla_{\V{v}} f &= 0\label{eq:VP1}, \\ \textbf{E} = -\nabla \varphi, \qquad -\Delta \varphi = q\int_{\mathbb{R}^3} f d\textbf{v}\equiv\rho\label{eq:VP2}.\end{align} \myblue{We note that in equation \ref{eq:VP2}, a background charge density $\rho_0$ is often subtracted in the context of periodic boundary conditions. In contrast, in this manuscript, we are interested in situations in which the charge density $\rho$ is the physical charge density, which may be the charge density of an isolated non-neutral beam. The only requirement on $\rho$ in this work is that it must have a compact support.}

\subsection{Particle-in-cell method}
In order to solve the system of equations \eqref{eq:VP1}--\eqref{eq:VP2}, PIC represents the distribution function as particles in a continuous phase space, each particle having a unique position $\textbf{X}_j$ and velocity $\textbf{V}_j$. The motion of the particles is described by Newton's law:\begin{equation}
    \dot{\textbf{X}}_j = \textbf{V}_j, \qquad \dot{\textbf{V}}_j = \frac{q}{m}(\textbf{E}_j+\textbf{V}_j\times \textbf{B}_0).
\end{equation} where $\textbf{E}_j$ indicates the self-consistent electric field $\textbf{E}$ evaluated at the particle position $\textbf{X}_j$. 

In the PIC method, $\textbf{E}$ is first computed on a grid, before it is evaluated at the particle positions via interpolation in order to advance the particles via Newton's law. Specifically, the charge density carried by the particles is spread onto grid points according to \begin{equation}
    \rho (\textbf{x}) = \sum_{j=1}^{N_p} qS(\textbf{X}_j-\textbf{x})\label{pic_rho}
\end{equation} where $S$ is the particle shape function. $S$ is usually chosen as a b-spline function of a certain order, which ensures exact charge conservation and retains positivity of the interpolation at any spatial location. Once the charge density $\rho$ is known at the grid points, one can compute $\textbf{E}$ and $\varphi$ by solving Poisson's equation on that grid using standard solvers. For example, for periodic boundary conditions and a uniform Cartesian grid, the electric field may be efficiently computed with Fast Fourier Transforms (FFTs):
\begin{align}
    \hat{\rho} (\textbf{k}) & = \sum_\textbf{x} \rho(\textbf{x}) \exp(-i\textbf{k}\cdot \textbf{x}),\label{e5} \\
    \hat{\varphi}(\textbf{k}) & = \frac{1}{k^2}\hat{\rho}(\textbf{k}), \qquad \varphi(\textbf{x}) = \sum_\textbf{k}\hat{\varphi}(\textbf{k}) \exp(i\textbf{k}\cdot \textbf{x}),\label{eq:spec_phi}\\
    \hat{\textbf{E}}(\textbf{k}) & = \frac{\textbf{k}}{ik^2}\hat{\rho}(\textbf{k}), \qquad \textbf{E}(\textbf{x}) = \sum_\textbf{k}\hat{\textbf{E}}(\textbf{k}) \exp(i\textbf{k}\cdot \textbf{x})\label{eq:spec_E}.
\end{align}


\myred{Here, we use the notation $\sum_\textbf{k}$ to denote the summation taken over a finite number of modes $\bold{k}$. Specifically, if we let $\bold{k}=[k_1, k_2]$, then \begin{equation}
    \sum_{\bold{k}}f(\bold{k}) = \sum_{k_1=-N_m/2}^{N_m/2-1}\sum_{k_2=-N_m/2}^{N_m/2-1} f([k_1, k_2])
\end{equation} for an arbitrary function $f(\bold{k})$, where $N_m$ is the number of modes in each dimension. This convention is followed throughout the paper.} Once the electric field is known at all grid points, it is evaluated at the particle positions by interpolation with the same shape function $S$: \begin{equation}
    \textbf{E}(\textbf{X}_j) = \sum_\textbf{x}\textbf{E}(\textbf{x})S(\textbf{x}-\textbf{X}_j).
\end{equation}

\subsection{Particle-in-Fourier scheme}
In order to improve energy conservation in the PIC method, a new scheme called particle-in-Fourier (PIF) has been proposed\cite{Evstatiev_2013, mitchell2019efficient}. In contrast to the standard PIC approach, the PIF scheme avoids grid heating effects by directly depositing the charge in Fourier space, as well as interpolating the electric field from the Fourier space to the particle locations. Let us assume our problem domain to be $\Omega=\bigotimes_d \left[-\frac{1}{2},\frac{1}{2}\right]$, where $d$ is the number of dimensions of the problem of interest, and denote $N_m$\footnote{For simplicity let us assume $N_m$ to be even.} as the number of Fourier modes in each dimension. If we perform a Fourier transform to \eqref{pic_rho}, we get 
\begin{align}
    \hat{\rho}(\textbf{k}) =& \int_\Omega q\sum_{j=1}^{N_p} S(\textbf{X}_j-\textbf{x})\exp(-i\textbf{k}\cdot \textbf{x})d\textbf{x},\\ =& q\hat{S}(\textbf{k})\sum_{j=1}^{N_p}\exp(-i\textbf{k}\cdot \textbf{X}_j), \label{e9}
\end{align}
where $\textbf{k}$ is of the form $2\pi\left(\sum_d n_d\hat{\textbf{d}}\right)$ for integers $n_d\in\left[-\frac{N_m}{2}, \frac{N_m}{2}-1\right]$, $\hat{\textbf{d}}$ being the unit vector in $d$-th dimension. Here, the Fourier coefficients of the shape function $\hat{S}(\textbf{k})$ can be evaluated analytically. Note that $\hat{\rho}(\textbf{k})$ in equation \eqref{e9} is different from the one in \eqref{e5}. This is because the b-spline functions are not band-limited in Fourier space, and in \eqref{e9} we have to truncate its analytical Fourier spectrum down to $N_m$ modes.  We then solve for the electric field using equation \eqref{eq:spec_E}. Once we obtain $\hat{\textbf{E}}(\textbf{k})$ and $\hat{\varphi}(\textbf{k})$ using the standard spectral method (equations \eqref{eq:spec_phi} and \eqref{eq:spec_E}), we evaluate the acceleration as \begin{equation}
    \textbf{a}_j = \frac{q}{m}\left(\textbf{V}_j \times \textbf{B}_0+\sum_{\textbf{k}} \hat{\textbf{E}}(\textbf{k})\hat{S}(\textbf{k})\exp(i\textbf{k}\cdot \textbf{X}_j)\right).\label{e10}
\end{equation} Note that both \eqref{e9} and \eqref{e10} require discrete Fourier transforms of data between equi-spaced points in Fourier space and non-uniformly distributed particle positions. This can be evaluated efficiently using the so-called non-uniform fast Fourier transform (NUFFT) \cite{Beylkin1995OnTF, DuttRokhlin}. Type 1 NUFFT evaluates the summation \begin{equation}
    \hat{f}(\textbf{k}) = \sum_{j=1}^{N_p}f(\textbf{X}_j)\exp(-i \textbf{X}_j\cdot\textbf{k}),
\end{equation} and its dual, Type 2 NUFFT, evaluates \begin{equation}
    f(\textbf{X}_j) = \sum_{\textbf{k}} \hat{f}(\textbf{k})\exp(i \textbf{X}_j\cdot \textbf{k}).
\end{equation} Both Type 1 and 2 NUFFT have a complexity of $\mathcal{O}(N_p+N_m^d\log N_m^d)$, which is much faster than a naive direct sum with complexity $\mathcal{O}(N_pN_m^d)$. The NUFFT thus makes the PIF scheme affordable to run long-time particle simulations with high numbers of particles and Fourier modes.

\subsection{Vico-Greengard-Ferrando free space Poisson solver}
\label{sec:vico}
Consider Poisson's equation \begin{equation}
    \Delta \varphi = -\rho,
\end{equation} in $\mathbb{R}^d$ with open (free space) boundary conditions, and where $\rho$ is a source distribution with a compact support $\Omega$. Without loss of generality, let us suppose that $\Omega=\bigotimes_d \left[-\frac{1}{2},\frac{1}{2}\right]$ is a $d$-dimensional unit box. The solution to Poisson's equation can then be written as \begin{equation}
    \varphi(\textbf{x}) = \int_\Omega g(\V{x}-\V{y})\rho(\textbf{y})d\textbf{y},\label{eq:green_convolution}
\end{equation}
where $g(\textbf{r})=\frac{1}{4\pi r}$ for $d=3$ and $g(\textbf{r})=\frac{-1}{2\pi}\log(r)$ for $d=2$. 

The slow decay of the Green's functions $g$ requires the application of fast algorithms for the efficient computation of the convolution in \eqref{eq:green_convolution}. For uniform Cartesian grids, this can be done with the convolution theorem: the Fourier transform of the convolution of the two functions in configuration space is the point-wise product of their Fourier transforms. Equation \eqref{eq:green_convolution} can therefore be evaluated efficiently by computing the Fourier transforms of $g$ and $\rho$ via FFTs, multiplying their Fourier transforms, and computing the inverse Fourier transform of the result, again with an FFT. This is what the Hockney-Eastwood scheme \cite{hockney1988computer,adelmann2019opal}, commonly used in PIC codes to deal with free space boundary conditions, does, by zero padding the charge distribution in each spatial direction on an interval of the same length as the original domain,
and defining a periodic Green’s function $\tilde{g}$ on that extended
domain. However, the scheme is only second-order accurate \cite{zou2021fftbased}, because of the second difficulty arising in \eqref{eq:green_convolution}, namely the fact that the Green's functions are singular at the origin. Hockney and Eastwood propose regularizing the Green's functions according to $\tilde{g}(\mathbf{0})=1$. The lack of smoothness of $\tilde{g}$ associated with this choice limits the convergence to second order \cite{rasmussen2011particle,hejlesen2019non}. Other elementary regularization procedures have been proposed \cite{budiardja2011parallel,chatelain2010fourier}, but they do not generally increase the order of convergence because they do not improve the smoothness of the regularized Green's function \cite{rasmussen2011particle,hejlesen2019non}.

The Vico-Greengard-Ferrando scheme also relies on a uniform grid, the convolution theorem and FFTs, but significantly improves the accuracy of the computation of the convolution in \eqref{eq:green_convolution} as compared to Hockney-Eastwood, by addressing the singularity of the Green's functions at the origin in a different manner. 

First, note that the Fourier transform of the Poisson Green's functions are in fact known analytically. However, a difficulty arises from the fact that in Fourier space, this transform diverges like $1/|\textbf{k}|^2$ at the origin. Vico, Greengard and Ferrando address that difficulty as follows. The key observation is that since we are only interested in the values of $\varphi$ and $\textbf{E}$ inside $\Omega$, one can subsitute the Green's functions $g$ with truncated Green's functions $g^{L}$ without affecting the results in $\Omega$. Specifically, let 
\begin{equation}
    g^L(\textbf{r}) = \begin{cases}
        \frac{1}{4\pi r}\textrm{rect}\left(\frac{r}{2L}\right)&\text{if } d=3\\ \\
        \frac{-1}{2\pi}\log(r)\textrm{rect}\left(\frac{r}{2L}\right)&\text{if } d=2
    \end{cases},
\end{equation}
where $\text{rect}(x)$ is the characteristic function on the unit interval: \begin{equation}
    \text{rect}(x)= \begin{cases}
        1&\text{for }|x|<1/2\\\\
        0&\text{otherwise}
    \end{cases}.
\end{equation}
If we set $L>\sqrt{d}$, then the solution $\varphi$ in $\Omega$ is exactly \begin{equation}
    \varphi(\textbf{x})=\int_\Omega g^L(\V{x}-\V{y})\rho(\textbf{y})d\textbf{y}.
\end{equation}
The advantage of truncating the Green's functions in physical space in this manner is that it removes the singularity of their Fourier transforms. Indeed, the Fourier transforms of the truncated Green's functions can also be computed analytically, and are the following smooth functions:\begin{equation}
    \mathcal{F}(g^L) = \hat{g}^L = \begin{cases}2\left(\frac{\sin(Ls/2)}{s}\right)^2& \text{if } d=3\\ \\ \frac{1-J_0(Ls)}{s^2}-\frac{L\log(L)J_1(Ls)}{s}& \text{if } d=2\end{cases}.
\end{equation} 
Because of the oscillatory nature of the truncated Green's functions $\hat{g}^L$, Vico, Greengard and Ferrando explain that there scheme in principle requires zero-padding to an extended domain $\tilde{\Omega}$ with a grid of size $(4N_g)^d$ in order to compute the convolution \eqref{eq:green_convolution} without aliasing errors. That would make their scheme more expensive than the Hockney-Eastwood method. However, they explain that the slightly more unfavourable scaling of the scheme may just be applied to a precomputation step to compute the numerical mollified Green's function $T$ on a grid of size $(2N_{g})^d$. Once computed, this numerical kernel $T$ can be used just like the periodic $\tilde{g}$ in the Hockney-Eastwood scheme, to compute $\varphi$ with the same computational complexity as in the Hockney-Eastwood scheme, via the convolution theorem and FFTs, and with much higher accuracy. Recently a modification taking advantage of the symmetry has been proposed in \cite{mayani2024massively}, which improves the memory storage costs of the Vico-Greengard-Ferrando scheme similar to that of the Hockney-Eastwood solver.

\section{Combining PIF with a spectral free space Poisson solver}
\label{sec:combine}
In this section, we present a novel way to combine PIF with the Poisson solver developed by Vico \textit{et al.}, \cite{Vico_2016} so that we retain the advantages of both, and can solve problems with free-space boundary conditions. Specifically, we can calculate the electric forces on particles efficiently and accurately for plasmas in free space, while achieving good conservation of charge, momentum, and energy.

\subsection{A gridless free space Poisson solver}
Since Vico \textit{et al.} solve Poisson's equation by convoluting the source with a mollified Green's function via the Fourier convolution theorem, we can view their solver as a spectral solver which is inherently compatible with PIF. Indeed, as shown in \eqref{e9}, in PIF we naturally  evaluate $\hat{\rho}=\mathcal{F}(\rho)$ required for the Vico-Greengard-Ferrando scheme with the NUFFT. Note that here, the charge density is again just
\begin{align}\label{eq:density_convol}
    \rho &= qS(\textbf{x})*\left(\sum_j \delta(\textbf{X}_j-\textbf{x})\right)
\end{align} 
where $(*)$ denotes the convolution operation. Since we now impose free space boundary condition for the Poisson solver, it is more natural to use radially symmetric shape functions. Hence, here $S(\textbf{x})$ is a radially symmetric function, as opposed to tensor product B-spline kernels commonly used in PIC \cite{birdsall2018plasma}. This is particularly convenient for the PIF scheme because one can analytically determine the Fourier coefficients of radially symmetric shape functions. \myblue{The use of radially symmetric shape functions has two main advantages, namely 1) they preserve rotational invariance, as it should be in the case of free space boundary conditions, and 2) the convolution between the mollified Green's function and the shape function is still radially symmetric, as we will see later in this section.} Some sample radially symmetric shape functions that we used in this work are provided in Appendix A.

For the solution to Poisson's equation, we may therefore write:
\begin{align}\label{eq:potential_convol}
    \varphi(\V{x}) &= g^L(\V{x}) * \rho(\V{x}) = q\left(g^L*S\right)(\V{x})*\left(\sum_j \delta(\textbf{X}_j-\textbf{x})\right),\\
    \V{E}(\V{x}) & = -\nabla \varphi(\V{x}).
\end{align}
Now, in order to determine the acceleration, we need to perform one more convolution between the fields and the shape function $S$. We may view the results as potential and electric force fields mollified by the shape function. For the sake of conciseness, let us denote these mollified fields as
\begin{align}
    \psi(\V{x}) &= (\varphi * S)(\V{x}),\\
    \mathcal{E}(\V{x}) & = (\V{E} * S)(\V{x}).
\end{align}
The accelerations are expressed as 
\begin{align}
    \V{a}_s & = \frac{q}{m}\left(\mathcal{E}(\V{X}_s) + \V{V}_s \times \textbf{B}_0\right),\\ 
    & = \frac{q}{m}\left(-(\nabla \varphi * S)(\textbf{X}_s) + \V{V}_s \times \textbf{B}_0\right),\\
    &= \frac{q}{m} \left[-q\nabla\left(g^L*S*S\right)*\left(\sum_j \delta(\textbf{X}_j-\textbf{x})\right) * \delta(\V{x}-\V{X}_s)+ \V{V}_s \times \textbf{B}_0\right].
\end{align}
For our scheme, we view $g^L * S$ and $\nabla\left(g^L*S*S\right)$ as our new convolution kernels in replacement of $g^L$, and apply the same Vico-Greengard-Ferrando scheme to compute the acceleration:  \begin{equation}
    \V{a}_s = \frac{q}{m}\left[-iq\sum_{\V{k}}\left(\V{k}\hat{g}^L(\V{k})\hat{S}^2(\V{k})\sum_j \exp(-i\textbf{k}\cdot \textbf{X}_j)\right)\exp(i\textbf{k}\cdot \textbf{X}_s)+\textbf{V}_s\times \textbf{B}_0\right],\label{eq:acc}
    \end{equation} 
    This can be performed efficiently using Type 1 \& 2 NUFFT. 
    
    We observe at this point that the explicit computation of the electric and potential fields are not required at any given time step of the PIF simulation, since Eq.\eqref{eq:acc} is sufficient to push the particles at any time step. The potential and electric force fields could therefore in principle only be computed as a post-processing step if desired, when analyzing the results of the simulations. For example, the potential field $\varphi$ can be computed using the following expression: \begin{equation}
    \hat{\varphi}(\textbf{k}) = q\hat{g}^L(\V{k})\hat{S}(\V{k})\sum_j \exp(-i\textbf{k}\cdot \textbf{X}_j)\label{eq:phi}.
\end{equation}

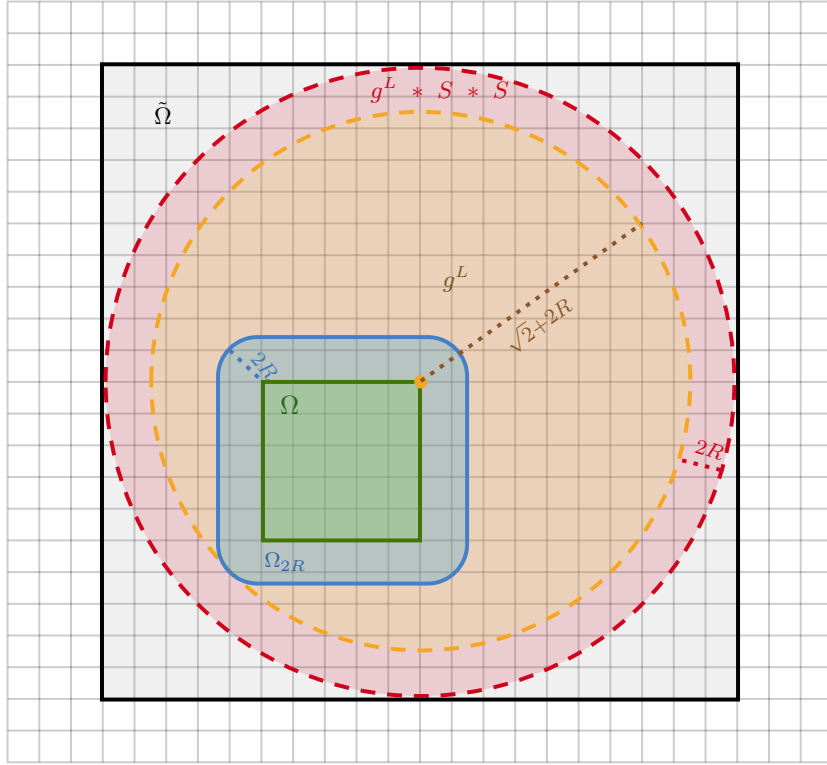
\begin{figure}
\centering
\tikzset{every picture/.style={line width=0.75pt}} 

\begin{tikzpicture}[x=0.6pt,y=0.6pt,yscale=-1,xscale=1]

\draw  [draw opacity=0] (90.52,60.18) -- (610.52,60.18) -- (610.52,540.18) -- (90.52,540.18) -- cycle ; \draw  [color={rgb, 255:red, 0; green, 0; blue, 0 }  ,draw opacity=0.2 ] (110.52,60.18) -- (110.52,540.18)(130.52,60.18) -- (130.52,540.18)(150.52,60.18) -- (150.52,540.18)(170.52,60.18) -- (170.52,540.18)(190.52,60.18) -- (190.52,540.18)(210.52,60.18) -- (210.52,540.18)(230.52,60.18) -- (230.52,540.18)(250.52,60.18) -- (250.52,540.18)(270.52,60.18) -- (270.52,540.18)(290.52,60.18) -- (290.52,540.18)(310.52,60.18) -- (310.52,540.18)(330.52,60.18) -- (330.52,540.18)(350.52,60.18) -- (350.52,540.18)(370.52,60.18) -- (370.52,540.18)(390.52,60.18) -- (390.52,540.18)(410.52,60.18) -- (410.52,540.18)(430.52,60.18) -- (430.52,540.18)(450.52,60.18) -- (450.52,540.18)(470.52,60.18) -- (470.52,540.18)(490.52,60.18) -- (490.52,540.18)(510.52,60.18) -- (510.52,540.18)(530.52,60.18) -- (530.52,540.18)(550.52,60.18) -- (550.52,540.18)(570.52,60.18) -- (570.52,540.18)(590.52,60.18) -- (590.52,540.18) ; \draw  [color={rgb, 255:red, 0; green, 0; blue, 0 }  ,draw opacity=0.2 ] (90.52,80.18) -- (610.52,80.18)(90.52,100.18) -- (610.52,100.18)(90.52,120.18) -- (610.52,120.18)(90.52,140.18) -- (610.52,140.18)(90.52,160.18) -- (610.52,160.18)(90.52,180.18) -- (610.52,180.18)(90.52,200.18) -- (610.52,200.18)(90.52,220.18) -- (610.52,220.18)(90.52,240.18) -- (610.52,240.18)(90.52,260.18) -- (610.52,260.18)(90.52,280.18) -- (610.52,280.18)(90.52,300.18) -- (610.52,300.18)(90.52,320.18) -- (610.52,320.18)(90.52,340.18) -- (610.52,340.18)(90.52,360.18) -- (610.52,360.18)(90.52,380.18) -- (610.52,380.18)(90.52,400.18) -- (610.52,400.18)(90.52,420.18) -- (610.52,420.18)(90.52,440.18) -- (610.52,440.18)(90.52,460.18) -- (610.52,460.18)(90.52,480.18) -- (610.52,480.18)(90.52,500.18) -- (610.52,500.18)(90.52,520.18) -- (610.52,520.18) ; \draw  [color={rgb, 255:red, 0; green, 0; blue, 0 }  ,draw opacity=0.2 ] (90.52,60.18) -- (610.52,60.18) -- (610.52,540.18) -- (90.52,540.18) -- cycle ;
\draw  [color={rgb, 255:red, 0; green, 0; blue, 0 }  ,draw opacity=1 ][fill={rgb, 255:red, 155; green, 155; blue, 155 }  ,fill opacity=0.15 ][line width=1.5]  (150.02,99.78) -- (551.02,99.78) -- (551.02,500.59) -- (150.02,500.59) -- cycle ;
\draw  [color={rgb, 255:red, 208; green, 2; blue, 27 }  ,draw opacity=1 ][fill={rgb, 255:red, 208; green, 2; blue, 27 }  ,fill opacity=0.15 ][dash pattern={on 5.63pt off 4.5pt}][line width=1.5]  (152.33,300.18) .. controls (152.33,190.73) and (241.07,102) .. (350.52,102) .. controls (459.98,102) and (548.71,190.73) .. (548.71,300.18) .. controls (548.71,409.64) and (459.98,498.37) .. (350.52,498.37) .. controls (241.07,498.37) and (152.33,409.64) .. (152.33,300.18) -- cycle ;
\draw  [color={rgb, 255:red, 245; green, 166; blue, 35 }  ,draw opacity=1 ][fill={rgb, 255:red, 237; green, 238; blue, 58 }  ,fill opacity=0.15 ][dash pattern={on 5.63pt off 4.5pt}][line width=1.5]  (181.12,299.68) .. controls (181.12,205.85) and (257.19,129.78) .. (351.02,129.78) .. controls (444.86,129.78) and (520.93,205.85) .. (520.93,299.68) .. controls (520.93,393.52) and (444.86,469.59) .. (351.02,469.59) .. controls (257.19,469.59) and (181.12,393.52) .. (181.12,299.68) -- cycle ;
\draw  [color={rgb, 255:red, 68; green, 128; blue, 200 }  ,draw opacity=1 ][fill={rgb, 255:red, 39; green, 160; blue, 218 }  ,fill opacity=0.26 ][line width=1.5]  (223.2,296.79) .. controls (223.2,283.04) and (234.35,271.9) .. (248.1,271.9) -- (355.32,271.9) .. controls (369.07,271.9) and (380.22,283.04) .. (380.22,296.79) -- (380.22,402.5) .. controls (380.22,416.25) and (369.07,427.4) .. (355.32,427.4) -- (248.1,427.4) .. controls (234.35,427.4) and (223.2,416.25) .. (223.2,402.5) -- cycle ;
\draw  [color={rgb, 255:red, 65; green, 117; blue, 5 }  ,draw opacity=1 ][fill={rgb, 255:red, 109; green, 204; blue, 63 }  ,fill opacity=0.25 ][line width=1.5]  (251.52,300.18) -- (350.52,300.18) -- (350.52,400.18) -- (251.52,400.18) -- cycle ;
\draw  [color={rgb, 255:red, 245; green, 166; blue, 35 }  ,draw opacity=1 ][fill={rgb, 255:red, 245; green, 166; blue, 35 }  ,fill opacity=1 ] (347.52,300.18) .. controls (347.52,298.25) and (349.09,296.68) .. (351.02,296.68) .. controls (352.96,296.68) and (354.52,298.25) .. (354.52,300.18) .. controls (354.52,302.12) and (352.96,303.68) .. (351.02,303.68) .. controls (349.09,303.68) and (347.52,302.12) .. (347.52,300.18) -- cycle ;
\draw [color={rgb, 255:red, 139; green, 87; blue, 42 }  ,draw opacity=1 ][line width=1.5]  [dash pattern={on 1.69pt off 2.76pt}]  (351.02,299.68) -- (490.52,200.18) ;
\draw [color={rgb, 255:red, 208; green, 2; blue, 27 }  ,draw opacity=1 ][line width=1.5]  [dash pattern={on 1.69pt off 2.76pt}]  (515.74,349.68) -- (541,356) ;
\draw [color={rgb, 255:red, 63; green, 122; blue, 190 }  ,draw opacity=1 ][line width=1.5]  [dash pattern={on 1.69pt off 2.76pt}]  (230.52,280.18) -- (251.52,300.18) ;

\draw (260.52,307.18) node [anchor=north west][inner sep=0.75pt]  [color={rgb, 255:red, 29; green, 107; blue, 19 }  ,opacity=1 ]  {$\si{\ohm}$};
\draw (250.52,406.18) node [anchor=north west][inner sep=0.75pt]  [font=\footnotesize,color={rgb, 255:red, 43; green, 99; blue, 167 }  ,opacity=1 ]  {$\si{\ohm}_{2R}$};
\draw (363,225) node [anchor=north west][inner sep=0.75pt]  [font=\small,color={rgb, 255:red, 111; green, 81; blue, 30 }  ,opacity=1 ]  {$g^{L}$};
\draw (317.52,106.18) node [anchor=north west][inner sep=0.75pt]  [font=\small,color={rgb, 255:red, 208; green, 2; blue, 27 }  ,opacity=1 ]  {$g^{L} \ *\ S\ *\ S$};
\draw (181,122) node [anchor=north west][inner sep=0.75pt]  [font=\small]  {$\tilde{\si{\ohm}}$};
\draw (399.06,272.17) node [anchor=north west][inner sep=0.75pt]  [font=\footnotesize,color={rgb, 255:red, 139; green, 87; blue, 42 }  ,opacity=1 ,rotate=-323.38]  {$\textcolor[rgb]{0.55,0.34,0.16}{\sqrt{2}}\textcolor[rgb]{0.55,0.34,0.16}{+2R}$};
\draw (524.92,333.18) node [anchor=north west][inner sep=0.75pt]  [font=\footnotesize,color={rgb, 255:red, 208; green, 2; blue, 27 }  ,opacity=1 ,rotate=-18.35]  {$2R$};
\draw (250.52,277.18) node [anchor=north west][inner sep=0.75pt]  [font=\footnotesize,color={rgb, 255:red, 62; green, 120; blue, 188 }  ,opacity=1 ,rotate=-50.33]  {$2R$};

\end{tikzpicture}
\caption{Our simulation domain based on the shape function and mollified Green's function chosen. Particles are within the domain $\Omega$ (green), and $S * S * \sum_{j}\delta(\textbf{x}-\textbf{X}_j)$ has a support $\Omega_{2R}$ (blue). \myred{The radius of the Green's function is chosen such that $g^L$ is able to cover the entire blue region $\Omega_{2R}$ whenever the center of $g^L$ is located inside the green region $\Omega$. To determine the smallest possible radius, we place the center of $g^L$ at the upper right corner of $\Omega$. After a simple calculation, we find the smallest radius to be $\sqrt{2}+2R$ for $g^L$.} Finally, the extended domain is $\tilde{\Omega}$ (gray).}
\label{fig:domain}
\end{figure}
The simplicity of the formulation described above partially conceals two minor difficulties, in that we need to be careful in choosing the size of the Fourier spectrum for the convolutions, and the radius of the mollified Green's function $g^L$. 

Let us start with the second point. Vico \textit{et al.} choose $L=1.5>\sqrt{d}$ in two dimensions so that the convolution between $\rho$ and $g^L$ gives the correct potential inside $\Omega$ (here we assume $d=2$ but the analysis applies to any dimension). In our solver, instead of $g^L$, the kernel function is $\nabla(g^L*S*S)$: our goal here is to correctly evaluate $\psi$ and $\mathcal{E}$ at any target position inside $\Omega$. Since the radius of the shape function is predetermined and the particle positions are given at each iteration, let us view $S * S * \sum_{j}\delta(\textbf{x}-\textbf{X}_j)$ as the source of the field. Its compact support $\Omega_{2R}$ is larger than the original domain $\Omega$ by an amount $2R$, as shown in figure \ref{fig:domain}. Following this observation, we choose the size of our mollified Green's function based on the following criterion: \begin{equation}
\forall \textbf{x}_0 \in \Omega, \quad \Omega_{2R} \subset D_{L}(\textbf{x}_0)
\end{equation} where $D_L (\textbf{x}_0)$ represents the disk of radius $L$ centered on $\textbf{x}_0$. In other words, when we translate the center of the mollified Green's function to any possible particle position inside $\Omega$, seen as a target, the distance from that target to any source in $\Omega_{2R}$ should be less than $L$, so that the mollified Green's function gives the same result as the true Green's function, and $\psi$ and $\mathcal{E}$ can be correctly evaluated at particle positions. It is easy to see that this is satisfied if $L \geq \sqrt{d} + 2R$, as shown in figure \ref{fig:domain}. When this condition is set for $L$, we note that the radius of $g^L*S$ is larger than $\sqrt{d}+3R$ and the radius of $g^L*S*S$ is larger than $\sqrt{d}+4R$. In practice, the radius of the shape function is small compared to $\sqrt{d}$, so that the modification to the minimal radius associated with $R$ is not significant. In our two-dimensional examples below, we chose $L=1.5$ as for the standard Poisson solver in two dimensions by Vico, Greengard and Ferrando.

Now, let us address the question of the required Fourier resolution. As discussed in Section \ref{sec:vico}, Vico \textit{et al.} explain that a zero-padding region extending to $\tilde{\Omega}=\bigotimes_{d}[-2,2]$, which is four times larger than the original domain $\Omega$, is necessary so that the aperiodic convolutions do not lead to aliasing, and so that the oscillatory Fourier
transforms of the convolution kernels are properly sampled. Our scheme has the same sampling requirements, unless we implement a precomputation step as suggested in \cite{Vico_2016} and discussed in the next section. To accurately compute our convolutions in the absence of the precomputation step, we thus upscale the resolution of the Fourier grid by a factor of 4. Although all our operations are done in Fourier space when we do not implement the precomputation step, in the remainder of this article we will sometimes refer to this upsampling in Fourier space as zero-padding of the domain $\Omega$ to the extended domain $\tilde{\Omega}=\bigotimes_{d}[-2,2]$, thereby using the same language as Vico \textit{et al.}. 

\subsection{On precomputation}
\label{sec:precomp}
As discussed in Section \ref{sec:vico}, in \cite{Vico_2016} Vico \textit{et al.} explain that in the common case in which the simulation domain and its grid remain fixed during the entire simulation, the extended computational grid with $(4N_g)^d$ grid points may only be used once, at the start of the simulation. This is to compute the real space representation $T$ on the computational grid of the convolutional kernel for the solution to Poisson's equation, which can be interpreted as a mollified Green's function associated with $g^L$. Once $T$ is computed, the electrostatic potential is calculated for any grid point in the domain $\Omega$ with the following convolution (for a two-dimensional setting)\begin{equation}\label{eq:mollified_convolution}
    \varphi(i,j) = \sum_{i,j}T(i-i', j-j')\rho(i,j)\;,
\end{equation}where $i,j$ are indices of the grid points. The advantage of precomputing $T$ in this manner is that Eq.\eqref{eq:mollified_convolution} may then be evaluated using standard aperiodic convolution, which only requires zero-padding by a factor of 2. As long as the simulation domain and the discretization remain the same, $T$ can be reused throughout the simulation, and only $(2N_g)^d$ grid points are needed for the calculation of the electrostatic potential.

Our PIF scheme does not lend itself naturally to the same precomputation method. Indeed, in our algorithm, the locations of the targets at which the potential must be evaluated are not known in advance, and change from one time step to the next. One therefore has two choices. The first choice is to rely on a quadruply extended grid in each direction as discussed above throughout the simulation, guaranteeing the absence of aliasing errors and spectral accuracy for our Poisson solver at each time step. The second choice is to precompute effective convolution operators $T_1$ and $T_2$ associated with the two kernels $g^L * S$ and $\nabla g^L * S * S$ via inverse FFTs for a quadruply refined Fourier grid in each direction, as in Vico \textit{et al.}, and then use $T_{1}$ and $T_{2}$ for standard aperiodic convolutions on a doubly extended grid in each direction for the entire simulation. This second approach is less expensive computationally. However, since $T_{1}$ and $T_{2}$ only lead to spectrally accurate quadrature rules for target locations matching those of the precomputation, this second approach introduces aliasing errors at the arbitrary target locations of our particles, as with standard FFT-based Poisson solvers \cite{hockney1988computer}. Accuracy is therefore reduced in the same way as for these solvers, with second order convergence to be expected. 
Both choices have merits, depending on the application, and on the number of Fourier modes and of particles for the simulation. In our numerical tests below, we will show results obtained for each of these two options. In most situations, the numerical error of particle based schemes such as ours is dominated by the noise associated with the Monte Carlo nature of the algorithm, and a second-order convergent Poisson solver provides sufficient accuracy. In such cases, relying on $T_{1}$ and $T_{2}$ and smaller padding regions does not lead to noticeable loss of accuracy, and is the most computationally efficient choice.

\myred{We note that the Hockney-Eastwood algorithm \cite{hockney1988computer} also has second order convergence, and for a standard PIC scheme has the advantage that an analytical expression for the Green's function is available in the physical space. One may then think that it could be considered as a simpler alternative to the Vico-Greengard solver with precomputation. However, for our PIF scheme, this advantage is irrelevant, as we have to convolve the Green's function with the shape function in Fourier space. In order to perform the convolution correctly, the usual extension of the Green's function to a $2L$ domain as in the standard Hockney-Eastwood solver will not be sufficient and the Green's function has to be zero padded adequately depending on the radius of the shape function before the convolution. This partially defeats the simplicity of the Hockney-Eastwood algorithm. On the other hand, the Vico-Greengard solver is a natural option for PIF schemes as the mollified Green's function is directly given as an analytical expression in Fourier space.}

\subsection{Dirichlet boundary conditions}\label{abc}
Following standard potential theory \cite{askham2017adaptive}, our free space particle-in-Fourier (FSPIF) scheme can be combined with Laplace solvers in order to simulate systems with Dirichlet boundary conditions. Consider the standard Poisson problem with Dirichlet boundary conditions:

\begin{equation}\label{eq:Dirichlet_goal}
    \begin{cases}
        \Delta \varphi = -\rho, & \text{in }\Omega_{\rho}\\\\
        \varphi =f &\text{on }\partial\Omega_{\rho}
    \end{cases},
\end{equation}
where $\partial\Omega_{\rho}$ is the boundary of $\Omega_{\rho}$. Let $\varphi^{P}$ be the free-space solution to the Poisson equation $\Delta\varphi^{P}=-\rho$ discussed previously, and let $\varphi:=\varphi^{P}+\varphi^{H}$, where $\varphi^{H}$ is a harmonic function satisying the following Laplace problem
\begin{equation}\label{eq:Laplace_BC}
    \begin{cases}
        \Delta \varphi^{H} = 0, & \text{in }\Omega_{\rho}\\\\
        \varphi^{H} =f-\varphi^{P} &\text{on }\partial\Omega_{\rho}
    \end{cases}.
\end{equation} Then $\varphi$ solves the desired Poisson problem with Dirichlet boundary conditions, equation \eqref{eq:Dirichlet_goal}. 

Here, we illustrate this generalization with a simple example, corresponding to two-dimensional plasmas confined within a disk on whose boundary one may apply a potential $f$, which may depend on the position along the disk boundary. For this particular situation, it is well-known that explicit solutions to equation \eqref{eq:Laplace_BC} may be constructed via elementary tools of complex analysis. Let us remind the reader of this construction. Suppose that $\varphi$ satisfies \begin{equation}
    \begin{cases}
        \Delta \varphi = -\rho, & \text{in }D_1(\textbf{0})\\\\
        \varphi =f &\text{on }C_1(\textbf{0})
    \end{cases},
\end{equation} 
As explained above, we split the solution $\varphi = \varphi^H+\varphi^P$, where $\varphi^P$ is the free space solution which we numerically solve according to \cite{Vico_2016}. The other component, $\varphi^H$, satisfies Laplace's equation on the unit disk:\begin{equation}
    \begin{cases}
        \Delta \varphi^H = 0, & \text{in } D_{1}(\textbf{0})\\\\
        \varphi^H = f -\varphi^P & \text{on }C_{1}(\textbf{0})
    \end{cases}.
\end{equation}
Poisson's formula then gives an explicit expression for $\varphi^H$ \cite{ahlfors1979complex}:\begin{equation}
    \varphi^H(\V{x}) = \frac{1}{2\pi}\int_{|\textbf{z}|=1}\frac{1-|\V{x}|^2}{|\textbf{z}-\V{x}|^2}(f-\varphi^P)(\V{z})d\theta. \label{eq:lap}
\end{equation} Writing $\V{z}=[u,v],\ \V{x}=[x,y]$,the components of the electric field at the particle position corresponding to $\varphi^H$ are
\begin{align}
    E^H_x &= \frac{1}{\pi}\int_{|z|=1}\frac{x[(u-x)^2+(v-y)^2]-[1-(x^2+y^2)](u-x)}{[(u-x)^2+(v-y)^2]^2}(f(\V{z})-\varphi^P(\V{z}))d\theta,\label{eq:EE1}\\
    E^H_y &= \frac{1}{\pi}\int_{|z|=1}\frac{y[(u-x)^2+(v-y)^2]-[1-(x^2+y^2)](v-y)}{[(u-x)^2+(v-y)^2]^2}(f(\V{z})-\varphi^P(\V{z}))d\theta.\label{eq:EE2}
\end{align} The integrals involved in the computation of the electrostatic potential and the electric field can easily be approximated numerically with the trapezoidal quadrature rule, which has spectral accuracy for analytic functions \cite{trefethen2014exponentially}. Let us assume that we put $N_B$ number of equi-spaced points $\textbf{z}_j$ on the boundary, we can then approximate $\varphi^H$ using the following equation:
\begin{equation}
    \varphi^H_{\text{num}}(\V{x}) = \frac{1}{N_B}\sum_{j=1}^{N_B}\frac{1-|\V{x}|^2}{|\textbf{z}_j-\V{x}|^2}(f-\varphi^P)(\V{z}_j), \label{eq:nlap}
\end{equation}
\begin{align}
    \notag
    E^H_{x,\text{num}} &= \frac{2}{N_B}\sum_{j=1}^{N_B}\frac{x[(u_j-x)^2+(v_j-y)^2]-[1-(x^2+y^2)](u_j-x)}{[(u_j-x)^2+(v_j-y)^2]^2}\\
    \label{eq:nEE1}
    &(f(\V{z}_j)-\varphi^P(\V{z}_j)),
\end{align}
\begin{align}
    \notag
    E^H_{y,\text{num}} &= \frac{2}{N_B}\sum_{j=1}^{N_B}\frac{y[(u_j-x)^2+(v_j-y)^2]-[1-(x^2+y^2)](v_j-y)}{[(u_j-x)^2+(v_j-y)^2]^2}\\
    \label{eq:nEE2}
    &(f(\V{z}_j)-\varphi^P(\V{z}_j)).
\end{align}
As long as the applied potential $f$ is smooth enough and the charge density $\rho$ is smooth enough for the free-space electrostatic potential $\varphi^P$ to be smooth, and as long as we do not evaluate the electrostatic potential too close to the boundary of the disk \cite{klockner2013quadrature}, we get high order convergence for $\varphi^{H}$ and $\V{E}$ with a simple application of the trapezoidal rule. For situations in which these circumstances are not satisfied, specialized quadrature rules could be implemented to maintain high order accuracy \cite{klockner2013quadrature}.

However, in order to determine the acceleration, we need to accurately evaluate $\psi^H = \varphi^H * S$ at particle positions. In our PIF framework, this is not trivial, because in order to compute convolutions in Fourier space, we need $\hat{\varphi}^H$, which is not readily available given the boundary integral form of $\varphi^{H}$. Here, we show that this convolution through Fourier space can be avoided if some additional conditions are met. To do so, we make the following assumptions: the shape function $S$ is spherically symmetric, and $\varphi$, which is originally defined on the unit disk, is smooth enough to be analytically extended to $D_{1+R}(\textbf{0})$ where $R$ is the radius of the shape function. Note that we do not have to evaluate this analytical extension, but we do need to assume its existence for the following derivation. Recall that $\psi = \varphi * S$, so that $\psi$ satisfies \begin{equation}\label{eq:Poisson_psi_noBC}
    \Delta \psi(\textbf{x}) = -(\rho * S)(\textbf{x}), \qquad \textbf{x}\in D_{1+R}(\textbf{0}).
\end{equation} In general,  it is a complicated matter to determine the boundary conditions that $\psi$ satisfy, and therefore solve Eq. \eqref{eq:Poisson_psi_noBC}. However, if we further assume that $\rho$ has a slightly smaller compact support $D_{1-R}(\textbf{0})$, which is often the case in applications \cite{davidson2001physics,Hurst2016}, then we have \begin{equation}
    \Delta \varphi = 0 \quad \text{in } \bigcup_{\V{x}\in C_1(\textbf{0})} D_R(\V{x}).
\end{equation} We can therefore use the mean value property of harmonic functions for $\varphi$ and the spherical symmetry of the shape function to obtain the following equalities: $\forall \mathbf{x}\in C_1(\V{0})$, \begin{align}
    \psi(\textbf{x}) &= \int_{D_R(\textbf{x})}S(\textbf{y}-\textbf{x})\varphi(\textbf{x})d\V{y},\\
    &=\int_0^R S(r) \int_{C_r(\textbf{x})} \varphi(\textbf{y})d\textbf{y}dr,\\
    &=\int_0^R 2\pi r\cdot S(r) \varphi(\V{x})dr,\\
    & = \varphi(\textbf{x}) = f.
\end{align} This identity provides the required boundary condition for $\psi(\textbf{x})$. We thus conclude that $\psi(\textbf{x})$ satisfies the following Poisson problem:
\begin{equation}\label{eq:Poisson_psi}
    \begin{cases}
        \Delta \psi(\V{x}) = -(\rho*S)(\V{x}), & \V{x}\in D_1(\textbf{0})\\\\
        \psi =f &\text{on }C_1(\textbf{0})
    \end{cases}.
\end{equation}
Let us define $\psi^P = \varphi^P * S$. Then we have \begin{equation}
    \Delta \psi^P = -\rho * S
\end{equation} with free space boundary conditions. If we call $\psi^H$ the solution to the problem \begin{equation}
    \begin{cases}
        \Delta \psi^H(\V{x}) = 0, & \V{x}\in D_{1}(\textbf{0})\\\\
        \psi^H = f-\psi^P &\text{on }C_{1}(\textbf{0})
    \end{cases}, \label{eq:pot}
\end{equation}
then $\psi=\psi^P+\psi^H$ satisfies Eq.\eqref{eq:Poisson_psi} as desired. 
Similarly, $\mathcal{E}=\textbf{E}*S$ can be found using equation \eqref{eq:EE1} and equation \eqref{eq:EE2} by replacing $f-\varphi^P$ with $f-\psi^P$. The potential energy due to $\varphi^H$ can also be simply determined numerically:
\begin{align}
    U_E^H &= \frac{1}{2}\int_\Omega \rho(\V{x}) \varphi^H(\V{x})d\V{x},\\
 &= \frac{q}{2}\sum_{\V{X}_j}\int_\Omega S(\V{X}_j-\V{x})\varphi^H(\V{x})d\V{x},\\
& = \frac{q}{2}\sum_{\V{X}_j}(\varphi^H*S)(\V{X}_j),\\
&=\frac{q}{2}\sum_{\V{X}_j}\psi^H(\V{X}_j).
\end{align} The simple methods of potential theory we present above do not only apply to the unit disk but to any geometry. For some two-dimensional domains, our approach can be combined with conformal mapping \cite{ahlfors1979complex,driscoll2002schwarz} for explicit or semi-explicit formulas. For other domains and for three-dimensional problems, standard numerical Laplace solvers may be used.\\

\subsection{Algorithm Outline}
We combine our gridless free-space Poisson solver with the standard PIF scheme to arrive at the following algorithm, outlined below
\begin{enumerate}
    \item Before the start of simulation, construct convolutional kernels $g^L*S$ and $\nabla (g^L*S*S)$ in Fourier space, and perform the precomputation step if desired;
    \item Initialize the particle positions, velocities and charges;
    \item Use Type 1 NUFFT to transform particle positions into an upscaled Fourier grid, in correspondence with the size of the extended domain $\tilde{\Omega}$. Let us relabel $\textbf{k} \rightarrow \textbf{k} / \alpha$ for the upscaled Fourier grid, then: \begin{equation}
        \hat{X}(\V{k}) = \sum_{j=1}^{N_p} \exp(-i\V{k}\cdot\V{X}_j);
    \end{equation}
    \item Find the corresponding Fourier modes of the free space potential $\varphi$ and acceleration $\V{a}$:
    \begin{align}
        \hat{\varphi}(\V{k}) &= q\hat{g}^L(\textbf{k})\hat{S}(\textbf{k})\hat{X}(\V{k}),\\
        \hat{\V{a}}(\V{k}) &= -\frac{i\textbf{k}}{m}q^2 \hat{g}^L(\V{k})\hat{S}^2(\textbf{k})\hat{X}(\V{k});
    \end{align}
    \item Find the free space acceleration of each particle using Type 2 NUFFT:
    \begin{equation}
        \V{a}_j^P = \sum_{\V{k}}\hat{\V{a}}(\V{k})\exp(i\V{k}\cdot\V{X}_j) + \frac{q}{m}\textbf{V}_j\times \textbf{B}_0.
    \end{equation} If free space boundary conditions apply to the electric field, then the true acceleration is $\V{a}_j=\V{a}_j^P$. Otherwise, if a Dirichlet boundary condition is imposed for the electric field, we need to additionally determine $\V{a}_j^H$ using a Laplace solver, and add $\V{a}_j^H$ to get the true acceleration $\V{a}_j=\V{a}_j^P+\V{a}_j^H$;
    \item Push the particles exactly as in standard PIC / PIF using the leapfrog algorithm if no magnetic field is applied, or using the Boris algorithm if there is a constant $\V{B}_0$.
\end{enumerate}
We make the following observations regarding the outline of the algorithm presented above. First, strictly speaking, the computation of $\hat{\varphi}$ in step 4. is not required to push the particles in step 6.  The computation of $\varphi$ and $\hat{\varphi}$ may therefore only be included at certain time steps if desired, or not at all, and the electrostatic potential computed as a post-processing step if needed to interpret simulations. Second, the leapfrog and Boris algorithms are suggested as suitable algorithms for particle pushing, but our scheme could be combined with other, possibly higher order, algorithms with suitable momentum and energy conservation properties.

\subsection{Energy and momentum conservation analysis}
\label{sec:error_analysis} 
\myblue{In this section, we provide energy and momentum conservation analyses for our free space particle-in-Fourier scheme, before and after time discretization.}
\subsubsection{\myblue{Energy conservation}}
\myblue{In \cite{mitchell2019efficient} the authors provided a proof for the energy conservation of the PIF scheme under the continuous time limit with periodic boundary conditions. Here we generalize this proof to free space boundary conditions. 
Without loss of generality, we let $q=m=1$. We define the potential energy as $U_E$, and the kinetic energy $U_K$, at time $t$ in terms of the Fourier modes of the relevant physical quantities as follows:\begin{align}
    U_E &= \frac{1}{2}\sum_\V{k} \overline{\hat{\rho}(\V{k})}\hat{\varphi}^P(\V{k}),\\
    U_K &= \frac{1}{2}\sum_j ||\V{V}_j||^2.\label{eq:E}
\end{align}
Taking time derivative of the first quantity $U_E$:
\begin{align}
    \frac{d}{dt} U_E &= \frac{d}{dt}\left(\frac{1}{2}\sum_{\bold{k}}\overline{\hat{\rho}(\bold{k})}\hat{\varphi}^P(\bold{k})\right),\\
    &=\frac{1}{2}\sum_{\bold{k}}G(\bold{k})\frac{d}{dt}\left(\hat{\rho}(\bold{k})\overline{\hat{\rho}(\bold{k})}\right),
\end{align} where $G$ is the Green's function we choose depending on the boundary condition: $G=\frac{1}{k^2}$ for periodic boundary conditions, and $G=\hat{g}^L$ for free space boundary conditions. Notice that since $\rho$ is a real function, $\hat{\rho}(\bold{k})=-\hat{\rho}(-\bold{k})$; Also, since the Green's function $g$ is symmetric by construction, $G(\bold{k}) = G(-\bold{k})$. Therefore, we have the following derivation:
\begin{align}
    \frac{d}{dt} U_E &=\frac{1}{2}\sum_{\bold{k}}G(\bold{k})\left(\hat{\rho}(\bold{k})\frac{d}{dt}\overline{\hat{\rho}(\bold{k})}+\overline{\hat{\rho}(\bold{k})}\frac{d}{dt}\hat{\rho}(\bold{k})\right),\\ &=\frac{1}{2}\sum_{\bold{k}}\left(G(\bold{k})\hat{\rho}(\bold{k})\frac{d}{dt}\overline{\hat{\rho}(\bold{k})}+G(-\bold{k})\overline{\hat{\rho}(-\bold{k})}\frac{d}{dt}\hat{\rho}(-\bold{k})\right),\\ &=\frac{1}{2}\sum_{\bold{k}}\left(G(\bold{k})\hat{\rho}(\bold{k})\frac{d}{dt}\overline{\hat{\rho}(\bold{k})}+G(\bold{k})\hat{\rho}(\bold{k})\frac{d}{dt}\overline{\hat{\rho}(\bold{k})}\right),\\ &=\sum_{\bold{k}}G(\bold{k})\hat{\rho}(\bold{k})\frac{d}{dt}\left(\overline{\hat{\rho}(\bold{k})}\right), \\ &=\sum_{\bold{k}}G(\bold{k})\hat{\rho}(\bold{k})\frac{d}{dt}\left(\overline{\sum_{j}\hat{S}(\bold{k})\exp(-i\bold{k}\cdot\bold{X}_j)}\right), \\ &= i\sum_{\bold{k}}\sum_j \bold{k}\cdot\bold{V}_j G(\bold{k})\overline{\hat{S}(\bold{k})}\hat{\rho}(\bold{k})\exp(i\bold{k}\cdot \bold{X}_j),\\
    \label{pe_der}
    &= i\sum_{\bold{k}}\sum_j \bold{k}\cdot\bold{V}_j G(\bold{k})\hat{S}(\bold{k})\hat{\rho}(\bold{k})\exp(i\bold{k}\cdot \bold{X}_j),
\end{align}
where the last line comes from the fact that $\hat{S}(\bold{k})$ is real.
Finally, we take the time derivative of the kinetic energy: \begin{align}
    \frac{d}{dt}U_K &= \frac{1}{2}\sum_{j}\frac{d}{dt}||\bold{V}_j||^2,\\
    &=\sum_j \bold{V}_j \cdot \frac{d\bold{V}_j}{dt},\\
    &=\sum_j \bold{V}_j\cdot\left(-S*\nabla\varphi(\bold{X}_j)\right),\\
    &= \sum_{j} \bold{V}_j \cdot\left(-\sum_{\bold{k}}\hat{S}(\bold{k})(i\bold{k}\hat{\varphi}(\bold{k}))\exp(i\bold{k}\cdot\bold{X}_j)\right),\\ 
    \label{ke_der}
    &=-\sum_{j}\sum_{\bold{k}}i\bold{k}\cdot\bold{V}_j \hat{S}(\bold{k})G(\bold{k})\hat{\rho}(\bold{k})\exp(i\bold{k}\cdot\bold{X}_j).
\end{align}
One can easily check the identity \begin{equation}
    \frac{d}{dt}(U_E+U_K)=0,
\end{equation}
from equations \eqref{pe_der} and \eqref{ke_der}.
Furthermore, we show that with the standard leapfrog time integrator and in the sole presence of electric forces, the energy conservation error has second order convergence for both periodic and free space boundary conditions. We begin by noticing that Eq.\eqref{eq:density_convol} and Eq.\eqref{eq:potential_convol} take the following form in Fourier space: }
\begin{equation}\hat{\rho}^n(\V{k})=\hat{S}(\V{k})\left(\sum_j \exp(-i\V{k}\cdot\V{X}_j^n)\right), \qquad \hat{\varphi}^n(\V{k})=G(\V{k})\hat{\rho}^n(\V{k}), \label{eq:id}
\end{equation} 
where the superscript $n$ denotes the quantities at time step $n$.
We therefore immediately find the potential energy to be:
\begin{equation}
    U_E^n = \frac{1}{2}\sum_{\V{k}}G(\V{k})\left|\hat{S}(\V{k})\right|^2\left|\sum_j \exp(-i\V{k}\cdot\V{X}_j^n)\right|^2.
    \label{eq:pe}
\end{equation}
Since the leapfrog scheme uses the particle velocities at half time steps, we write:
\begin{align}
    \V{V}^{n+1/2}_s &= \V{V}_s^{n-1/2}+\Delta t \cdot\V{a}_{s}^n,\\
    &= \V{V}_s^{n-1/2}+\Delta t\sum_\V{k}\left(-i\V{k}G(\V{k})\left|\hat{S}(\V{k})\right|^2\exp(i\V{k}\cdot\V{X}^n_s)\sum_j\exp(-i\V{k}\cdot\V{X}_j^n)\right), \label{eq:velocityupdate}
\end{align} where the expression for $\V{a}_s^n$ was introduced in equation \eqref{eq:acc}. If we use the midpoint rule to define $\V{V}^n$, then we have the following relationship:\begin{equation}
    \V{V}_s^{n+1/2} = \V{V}_s^{n}+\frac{\Delta t}{2}\sum_\V{k}\left(-i\V{k}G(\V{k})\left|\hat{S}(\V{k})\right|^2\exp(i\V{k}\cdot\V{X}^n_s)\sum_j\exp(-i\V{k}\cdot\V{X}_j^n)\right).
\end{equation} We can then compute particle positions at time $n+1$ in terms of positions and velocities at time $n$:
\begin{align}
    \V{X}_s^{n+1} &= \V{X}_s^n +\Delta t \V{V}_s^{n+1/2},\\
    &=\V{X}_s^n +\Delta t \V{V}_s^n +\frac{\Delta t^2}{2}\sum_{\V{k}}\left(-i\V{k}G(\V{k})\left|\hat{S}(\V{k})\right|^2\exp(i\V{k}\cdot\V{X}^n_s)\sum_j\exp(-i\V{k}\cdot\V{X}_j^n)\right),\\
    &=\V{X}_s^n +\Delta t \V{V}_s^n +\mathcal{O}(\Delta t^2).
\end{align}
We can plug this expression into the equation \eqref{eq:pe} for the potential energy:
\begin{align}
    U_E^{n+1} &= \frac{1}{2}\sum_{\V{k}}G(\V{k})\left|\hat{S}(\V{k})\right|^2\left|\sum_s \exp\left(-i\V{k}\cdot\left(\V{X}_s^n +\Delta t \V{V}_s^n +\mathcal{O}(\Delta t^2)\right)\right)\right|^2. \label{exp1}
\end{align}
Now, for $\Delta t$ small, 
\begin{align}
    \exp\left(-i\V{k}\cdot\left(\V{X}_s^n +\Delta t \V{V}_s^n +\mathcal{O}(\Delta t^2)\right)\right) &=\exp(-i\V{k}\cdot\V{X}_s^n)(1-i\Delta t\V{k}\cdot \V{V}_s^n+\mathcal{O}(\Delta t^2)).
\end{align}
Substituting this equality in equation \eqref{exp1} gives
\begin{align}
     U_E^{n+1} &= \frac{1}{2}\sum_{\V{k}}G(\V{k})\left|\hat{S}(\V{k})\right|^2\left|\sum_s \exp\left(-i\V{k}\cdot\V{X}_s^n\right)\left(1-i\Delta t\V{k}\cdot\V{V}_s^n+\mathcal{O}(\Delta t^2)\right)\right|^2,\\
     \notag
    &=U_E^n - \Delta t \sum_{\V{k}}G(\V{k})\left|\hat{S}(\V{k})\right|^2 \overline{\left(\sum_s i\V{k}\cdot\V{V}_s^n \exp(-i\V{k}\cdot\V{X}_s^n)\right)}\\
    &\left(\sum_j \exp (-i\V{k}\cdot \V{X}_j^n)\right)+ \mathcal{O}(\Delta t^2),\\
    \notag
    & = U_E^n - \sum_s\V{V}_s^{n}\cdot\Delta t\sum_\V{k}\left(-i\V{k}G(\V{k})\left|\hat{S}(\V{k})\right|^2\sum_j \exp(i\V{k}\cdot(\V{X}_s^n-\V{X}_j^{n}))\right)\\
    \label{eq:pe_final}
    &+\mathcal{O}(\Delta t^2).
\end{align}
The velocities at the next time step, $\V{V}_s^{n+1}$, are calculated as:

\begin{align}
    \V{V}^{n+1}_s& = \frac{1}{2}(\V{V}^{n+1/2}_s + \V{V}^{n+3/2}_s), \\ 
    &=\frac{\V{V}_s^{n+1/2}}{2}+\frac{1}{2}(\V{V}_s^{n+1/2} + \Delta t \cdot \V{a}_s^{n+1}),\\
    &=\V{V}_s^{n+1/2}+\frac{\Delta t}{2} \V{a}_s^{n+1},\\
    & = \V{V}_s^{n+1/2}+\frac{\Delta t}{2} \sum_\V{k}\left(-i\V{k}G(\V{k})\left|\hat{S}(\V{k})\right|^2\exp(i\V{k}\cdot\V{X}_s^{n+1})\sum_j \exp(i\V{k}\cdot-\V{X}_j^{n+1})\right),\\
    & = \V{V}_s^{n+1/2} +\frac{\Delta t}{2}\\
    &\sum_\V{k}\left(-i\V{k}G(\V{k})\left|\hat{S}(\V{k})\right|^2\sum_j \exp\left(i\V{k}\cdot(\V{X}_s^n-\V{X}_j^n)+i\Delta t\V{k}\cdot(\V{V}_s^{n+1/2}-\V{V}_j^{n+1/2}) \right)\right),\\
    &=\V{V}_s^{n+1/2}+\frac{\Delta t}{2} \sum_\V{k}\left(-i\V{k}G(\V{k})\left|\hat{S}(\V{k})\right|^2\sum_j \exp\left(i\V{k}\cdot(\V{X}_s^n-\V{X}_j^{n})\right)\right) +\mathcal{O}(\Delta t^2),\\
    & = \V{V}_s^{n}+\Delta t\sum_\V{k}\left(-i\V{k}G(\V{k})\left|\hat{S}(\V{k})\right|^2\sum_j \exp(i\V{k}\cdot(\V{X}_s^n-\V{X}_j^{n}))\right) +\mathcal{O}(\Delta t^2).
\end{align}
Note that we again used the midpoint rule for the last equality. This equality gives us the total kinetic energy at time $n+1$ to second order in $\Delta t$:
\begin{align}
    U_K^{n+1} &= \frac{1}{2}\sum_s ||\V{V}_s^{n+1}||^2,\\
    \notag
    &=U^n_K+\Delta t\sum_s \V{V}_s^n \cdot\sum_{\V{k}}\left(-i\V{k}G(\V{k})\left|\hat{S}(\V{k})\right|^2\sum_j \exp(i\V{k}\cdot(\V{X}_s^n-\V{X}_j^n))\right)\\
    \label{eq:ke_final}
    &+\mathcal{O}(\Delta t^2).
\end{align}
Adding equations \eqref{eq:ke_final} and \eqref{eq:pe_final} we conclude that we have second order convergence for energy conservation with respect to $\Delta t$:\begin{equation}U_K^{n+1}+U_E^{n+1}=U_K^n+U_E^n+\mathcal{O}(\Delta t^2).
\end{equation}
\myred{This result, showing second order convergence for energy conservation independently of the number of particles and modes, is not shared by classical particle-in-cell scheme. In a standard explicit particle-in-cell scheme, the grid heating effect and finite grid instability will cause energy to increase without bound over time if the spatial grid does not resolve the Debye length.} We note that for our numerical tests, we will consider physical situations that involve other fields, such as an external magnetic field $\V{B}_0$ or a homogeneous potential field $\varphi^H$ derived from Dirichlet boundary conditions. By choosing suitable numerical algorithms, we can preserve the second order convergence of energy conservation in these situations. In the presence of a magnetic field $\V{B}_0$, one may choose the Boris algorithm to update the particle positions and velocities, which achieves exact energy conservation if $\V{E}=\V{0}$ and at least second order otherwise \cite{hairer2018energy}. In the presence of a homogeneous potential field $\varphi^H$ to satisfy Dirichlet boundary conditions, one can use the boundary integral method mentioned in Section \ref{abc} to compute the full electric field. The mathematical proof of energy conservation for these situations is left as future work; the numerical convergence of energy conservation in these cases will however be confirmed in the next section.
\subsubsection{\myred{Momentum conservation}}
\label{theo_mom_cons}
\myred{In this section, we analyze the momentum conservation of the free space particle-in-Fourier scheme. Without loss of generality we assume the mass of a particle is $m=1$. We define the total momentum of the system at time $t^{n-1/2}$ as \begin{equation}
    \bold{p}^{n-1/2}=\sum_{i=1}^{N_p} \bold{V}_i^{n-1/2}.
\end{equation}
By equation \ref{eq:velocityupdate}, we are able to determine the momentum at the next time step: \begin{equation}
    \bold{p}^{n+1/2} = \bold{p}^{n-1/2} + \Delta t\sum_{j=1}^{N_p}\exp(-i\bold{k}\cdot \bold{X}_j^n)\left(\sum_{\bold{k}}\left(-i\bold{k}G(\bold{k})|\hat{S}(\bold{k})|^2\sum_{i=1}^{N_p}\exp(i\bold{k}\cdot \bold{X}_i^n)\right)\right).
\end{equation} Therefore, the difference between the momentums at two adjacent time steps are 
\begin{align}
    \notag
    \bold{p}^{n+1/2}-\bold{p}^{n-1/2}&=\Delta t\sum_{j=1}^{N_p}\exp(-i\bold{k}\cdot \bold{X}_j^n)\\
    &\left(\sum_{\bold{k}}\left(-i\bold{k}G(\bold{k})|\hat{S}(\bold{k})|^2\sum_{i=1}^{N_p}\exp(i\bold{k}\cdot \bold{X}_i^n)\right)\right),\\
    &=-i\Delta t\sum_{\bold{k}}\bold{k}G(\bold{k})|\hat{S}(\bold{k})|^2\sum_{i=1}^{N_p}\sum_{j=1}^{N_p} \exp(i\bold{k}\cdot(\bold{X}_i^n-\bold{X}_j^n)). \label{eq:momentumsum}
\end{align}
Finally, we make note that $G(\bold{k})=G(-\bold{k})$, $\hat{S}(\bold{k})=\hat{S}(-\bold{k})$, and \begin{align}
    \mathcal{B}(\bold{k})&\equiv\sum_{i=1}^{N_p}\sum_{j=1}^{N_p}\exp \left(i\bold{k}\cdot(\bold{X}_i^n-\bold{X}_j^n)\right),\\
    &=\sum_{j=1}^{N_p}\sum_{i=1}^{N_p}\exp \left(i\bold{k}\cdot(\bold{X}_i^n-\bold{X}_j^n)\right),\\
    &=\sum_{j=1}^{N_p}\sum_{i=1}^{N_p}\exp \left(-i\bold{k}\cdot(\bold{X}_j^n-\bold{X}_i^n)\right),\\
    &=\sum_{i=1}^{N_p}\sum_{j=1}^{N_p}\exp \left(-i\bold{k}\cdot(\bold{X}_i^n-\bold{X}_j^n)\right),\\
    &=\mathcal{B}(-\bold{k}),
\end{align} are all even functions of $\bold{k}$, whereas $\bold{k}=-(-\bold{k})$ is an odd function of $\bold{k}$. We therefore conclude that the summation in equation \ref{eq:momentumsum} equals zero, and thus the total momentum of the system is conserved over time.
}
\subsection{\myblue{Comparison between PIC and PIF schemes}}

\myblue{The computational cost of PIF schemes scales as $\mathcal{O}((|log\varepsilon|+1)^dN_p + N_m^dlogN_m^d)$, where $\varepsilon$ is the tolerance used for the NUFFT. In contrast, the computational cost of a standard PIC scheme with cloud-in-cell shape function and an FFT-based field solver scales as $\mathcal{O}(2^dN_p + N_m^dlogN_m^d)$. Thus the standard PIC scheme has a computational complexity similar to the coarsest NUFFT PIF scheme with the tolerance of $\varepsilon=10^{-1}$. In other words, for the same total number of particles and number of modes, PIF schemes are costlier than the standard PIC scheme. However, the PIF scheme has advantages in the following respects:
\begin{enumerate}
\item As mentioned in \cite{mitchell2019efficient}, we can use pre-computed shape functions of arbitrary order in PIF schemes without any increase in their computational cost, whereas for PIC schemes the computational cost increases with the order of the shape function.
\item We can choose the tolerance $\varepsilon$ of the NUFFT in a PIF scheme depending on the desired levels of conservation in energy and charge for the problem of interest. In practice, even with the same number of particles and grid points/modes, PIF schemes can only be slightly costlier than the PIC schemes as shown in Section \ref{grid_heating}. 
\item The mesh size in explicit PIC schemes have to be on the order of the Debye length in order to avoid grid heating. This typically requires a much higher number of grid points than the number of modes considered in PIF schemes. Moreover, in order to have the same number of particles per cell, the total number of particles have to be increased as well. As a result of these practical constraints, one ends up having a much higher number of particles and grid points in an explicit PIC scheme compared to the PIF scheme, which in turn leads to increased time to solution. As an example we shown in Section \ref{grid_heating} that for comparable accuracies in the quantities of interest, the PIF scheme leads to more than an order of magnitude speedup compared to the explicit PIC scheme we implemented for comparison.
\item PIF schemes are spectral particle methods, which belong to the general class of geometric, structure preserving PIC schemes \cite{Evstatiev_2013}. There are other geometric PIC schemes \cite{campos2022variational,squire2012geometric,kraus2017gempic,he2016hamiltonian} which also offer excellent long time conservation and stability properties. Compared to them, PIF schemes retain the intuitive and ease of implementation nature of standard explicit PIC schemes. A more detailed comparison between PIF schemes and other geomteric PIC schemes in terms of accuracy and performance is beyond the scope of current work and will be carried out elsewhere in the future.   

\end{enumerate}

In terms of parallelization strategies, usually domain decomposition is used for PIC schemes where both the particles and grid points are divided between the cores. PIF schemes are relatively more global in nature than PIC schemes. Hence, as mentioned in \cite{mitchell2019efficient} they are more suited for particle decomposition parallelization strategy where only the particles are divided between cores and each core carries all the Fourier modes. However, this strategy has a bottleneck in both memory and computing costs when the simulation requires a lot of Fourier modes (as they are duplicated and not decomposed).  In order to alleviate this issue recently a space-time parallelization strategy has been proposed in \cite{muralikrishnan2024parapif} which shows promise for large scale 3D-3V PIF simulations on extreme-scale computing architectures.}

\section{Numerical Tests}\label{sec:tests}
In this section, we first test the numerical accuracy of our gridless free space Poisson solver with the method of manufactured solutions. We then verify our free space PIF implementation in a two-dimensional beam dynamics problem which has been extensively studied analytically \cite{cerfon2013analytic,Cerfon_2016}. Finally, the accuracy of our PIF scheme for Dirichlet boundary conditions based on FSPIF and the boundary integral method is also examined. For NUFFTs needed for these examples, we rely on \texttt{fiNUFFT}, an efficient library presented in \cite{barnett2019parallel}. The source codes for all our numerical results can be found at \url{https://github.com/Nigel-Shen/Free-Space-Particle-In-Fourier}.
\subsection{Free space Poisson problem}\label{fspp}
In this test case, we place 30 random particles in the unit box, and solve the free space Poisson problem using our gridless solver. We choose the shape function to be a truncated Gaussian with standard deviation $\sigma=1/100$ and truncation radius $R=1/8$ (see Appendix A), so that $\rho$ is almost $C^{\infty}$ within the domain $\Omega$. We choose the radius of the mollified Green's function to be $L=1.75$, and let the extended box, which includes the zero padded regions, be $\tilde{\Omega}=[-2,2]\times[-2,2]$. \\

\begin{figure}
	\centering
    \includegraphics[width=0.6\textwidth]{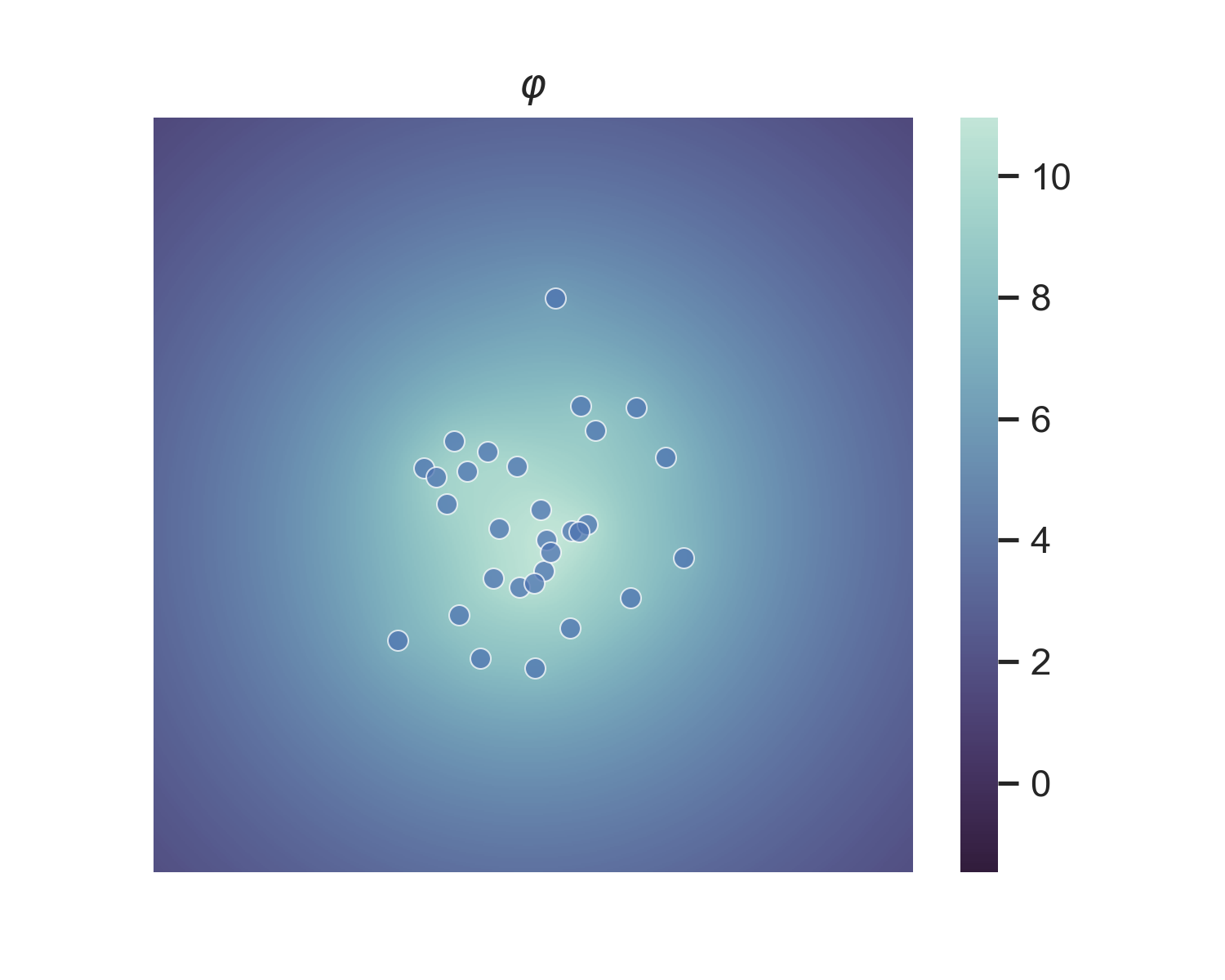}
    \hfill
	\includegraphics[width=1\textwidth]{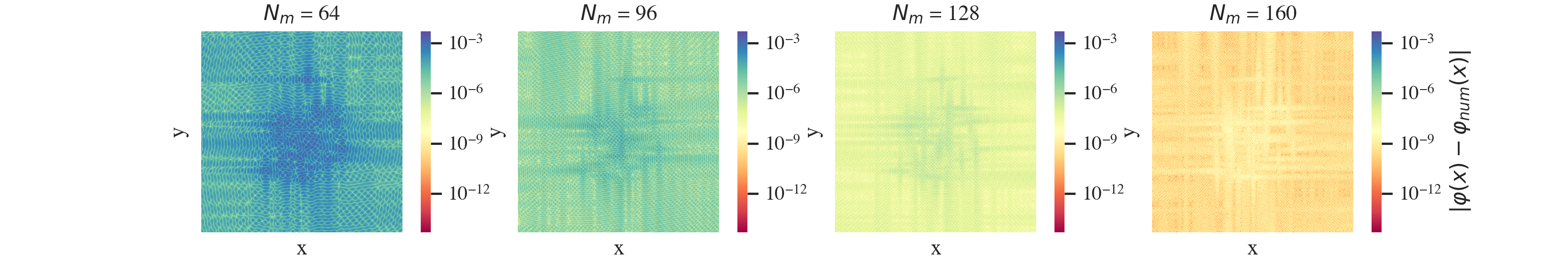}
  \caption{Potential $\varphi$ from the analytic solution given by equation \eqref{eq:ana} with particles represented as blue circles (top) and numerical error of our free space Poisson solver using different numbers of modes (bottom). We choose the radius of the Green's function to be $L=1.75$, and the solver uses a $(4N_m)^2$ box.}
	\label{fig:fig3}
\end{figure}
\begin{figure}
	\centering
        \includegraphics[width=0.8\textwidth]{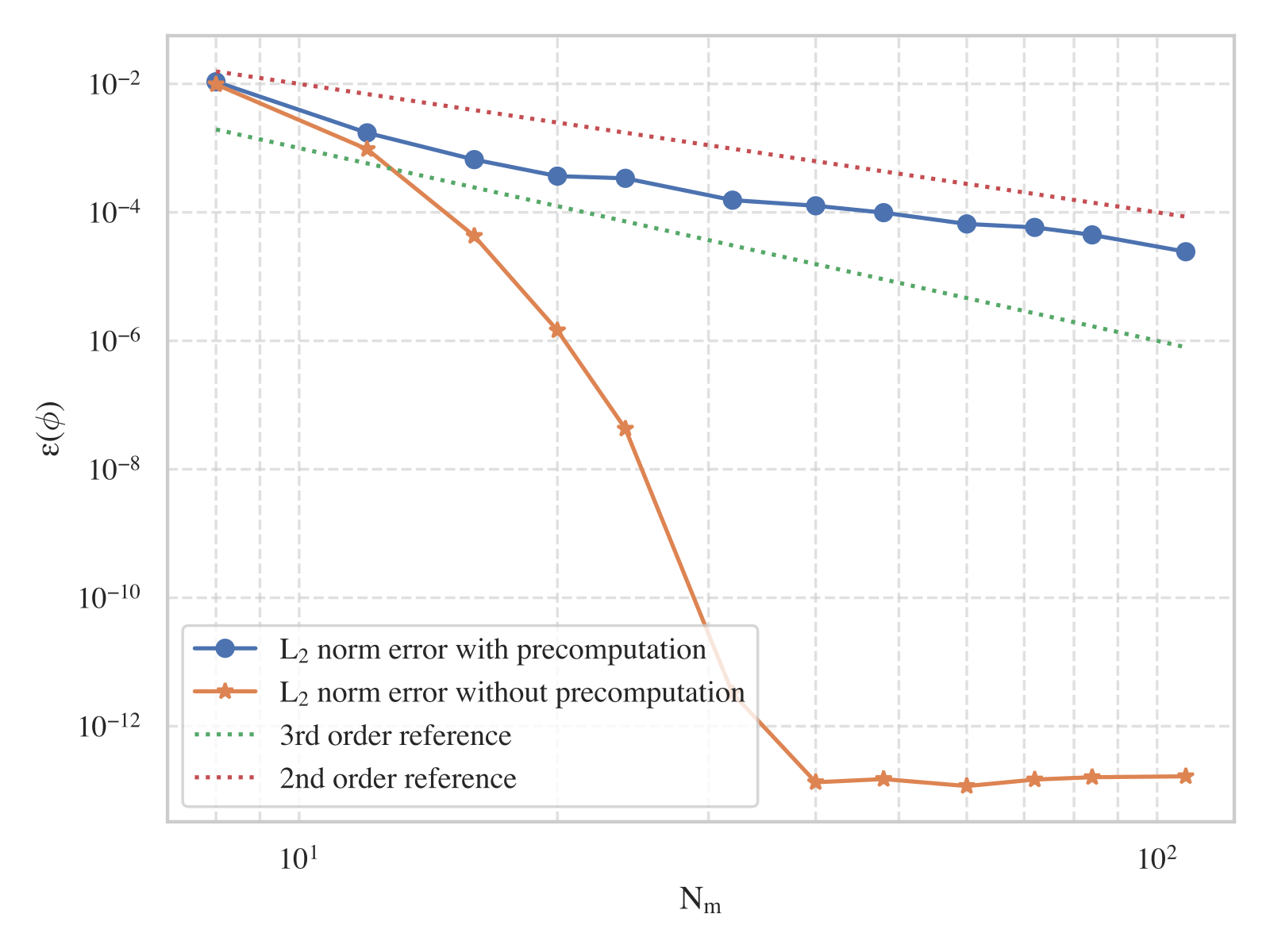}
	\caption{Global absolute errors in $L_2$ norms versus the number of modes in each dimension $N_m$, calculated in an upscaled physical grid of size $512\times 512$ using Fourier interpolation. Notice that without performing the precomputation step, the solver achieves spectral accuracy, whereas with the precomputation step, it has second order convergence.}
        \label{fig:fig4}
\end{figure}

Let us denote the particle positions as $\textbf{X}_j$. Both $\rho$ and $\varphi$ can be determined analytically:
\begin{align}
    \rho(\textbf{x}) = q\sum_j S(\textbf{X}_j-\textbf{x}) = \frac{q}{2\pi\sigma^2} \sum_j \exp\left(-\frac{||\textbf{X}_j-\textbf{x}||^2}{2\sigma^2}\right),\\
    \varphi(\textbf{x}) = \frac{q}{4\pi}\sum_j-\left(\text{Ei}\left(\frac{||\textbf{X}_j-\textbf{x}||^2}{2\sigma^2}\right)-\log(||\textbf{X}_j-\textbf{x}||^2)\right)\label{eq:ana}
\end{align} \myred{where Ei is the exponential integral function defined as \begin{equation}
    \text{Ei}(x)=-\int_{-x}^\infty\frac{e^{-t}}{t}dt.
\end{equation}}
We vary the number of modes $N_m$ and evaluate the numerical solution $\varphi(\V{x})$ on a refined grid of resolution $128\times128$ using Fourier interpolation (\textit{i.e.}, padding zeros in the spectral space before inverse FFT), and compare it with the analytical solution.\\

We begin by plotting the analytical solution $\varphi$ and compare it with the numerical solutions obtained when $N_m= 16, 24, 32$ and $40$ in Figure \ref{fig:fig3}. We then plot the absolute error $||\varepsilon||_2 = ||\varphi-\varphi_{num}||_2$ as a function of the number of modes $N_m$ on a log-log scale, shown in Figure \ref{fig:fig4}. Spectral accuracy can be clearly observed as the $L_2$ norm decays faster than any polynomial order in the case without precomputation step. If we do a precomputation step before solving the system, the solver now is second order if evaluated at arbitrary positions inside the domain, for the reason explained in Section \ref{sec:precomp}. 

\subsection{Infinitely long beam confined by a guiding magnetic field, in free space} 
\label{cb}
In this example, we consider an infinitely long charged particle beam in the $z$ direction, in free space, and confined by a strong magnetic field aligned with the $z$ axis. This situation has been considered as an insightful model for understanding the space charge dynamics in high intensity cyclotrons, and in particular the process of beam axisymmetrization \cite{yang2010beam,cerfon2013analytic,Cerfon_2016}. We let our domain be $\Omega=[-1/2, 1/2]\times [-1/2, 1/2]$. The initial extended box $\tilde{\Omega}$ including the zero-padded regions is set to be $[-2,2]\times[-2,2]$. A precomputation step is optionally performed, which then enables us to only consider $[-1,1]\times[-1,1]$ as the computational box including the zero padded regions. We normalize the spatial length unit by the Debye length $\lambda_D=\sqrt{\varepsilon_0 T/q\rho}$, and the temporal unit by the inverse plasma frequency $\omega_p^{-1}=\sqrt{\varepsilon_0 m/q\rho}$, where the particle mass $m$ and charge $q$ are each normalized to one. In order to confine our charged particles, we apply an external constant magnetic field $\textbf{B}_0=300\mathbf{z}$. The magnitude of $\textbf{B}_0$ has been chosen in such a way that the force due to the magnetic field is approximately 50 to 100 times larger in magnitude compared to the force due to the electric field. The initial particle distribution is an elongated normal distribution in configuration space, and Maxwellian in velocity space:
\begin{equation}
    f(x,y,v_x,v_y,0) = \frac{n(x,y)}{2\pi}\exp\left(-\frac{v_x^2+v_y^2}{2}\right)
\end{equation} where \begin{equation}
    \label{eq:ic_cyclotron}
    n(x,y)=\frac{1}{\pi\sigma_x\sigma_y}\exp\left(-\frac{x^2}{\sigma_x^2}\right)\exp\left(-\frac{y^2}{\sigma_y^2}\right) \text{ with }\sigma_x=\frac{1}{30}, \sigma_y=\frac{1}{10}.
\end{equation}

\begin{figure}
	\centering
	\includegraphics[width=1\textwidth]{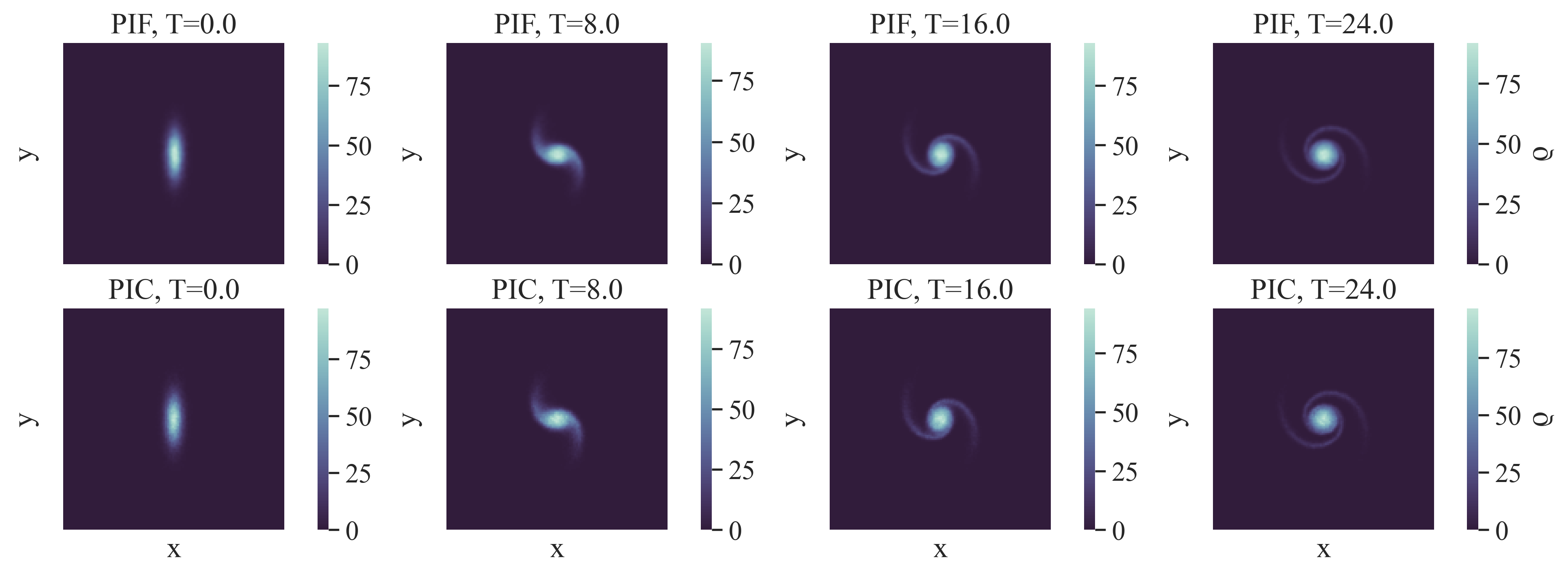}
	\caption{Charge density $\rho$ at different times of an infinitely long non-neutral beam confined by a strong constant and uniform magnetic field, simulated using 40,000 particles and $N_m=N_g=128$ using free space PIC and PIF methods. For the FSPIF scheme, the precomputation step is not performed.}
	\label{fig:fig5}
\end{figure}

\begin{figure}
	\centering
	\includegraphics[width=1\textwidth]{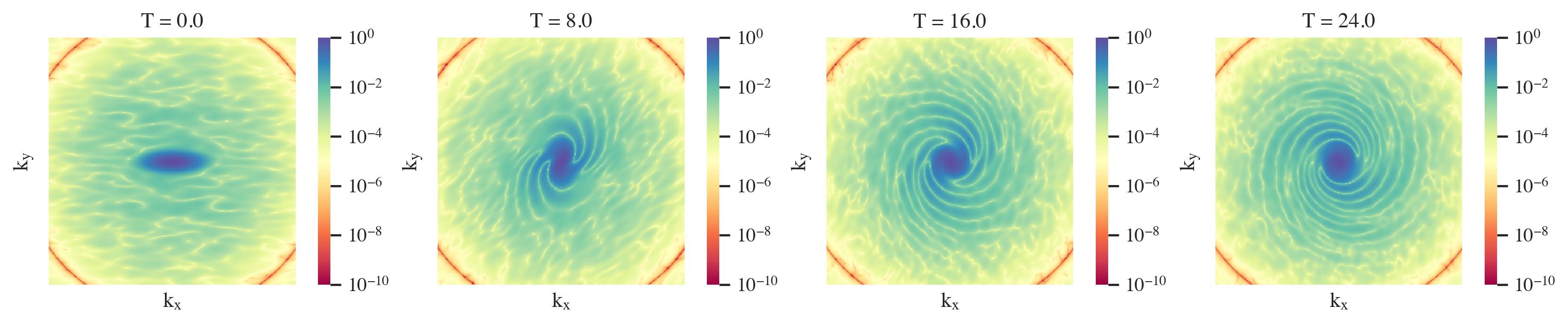}
	\caption{Absolute values of Fourier modes when $N_m=128$ for the same nonneutral beam simulation as in Figure \ref{fig:fig5}. These figures indicate that most high-frequency modes are not very significant, and there is no obvious mode-coupling or aliasing effects.}
	\label{fig:fig6}
\end{figure}

We fix the number of particles to be 40,000, and use a second-order radially symmetric b-spline function for the shape function (see Appendix A). The radius of the shape function is set to be $R=1/N_m$ where $N_m$ is the Fourier resolution. We first perform the simulation with $N_m=128$ in FSPIF and compare it with the result obtained using the standard PIC method. For FSPIC, the charges are spread onto the grid points using the second order tensor-product b-spline function, the standard Vico solver \cite{Vico_2016} is applied for spectrally accurate field solutions at the grid points, and the force is interpolated back to the particles using the same shape function as the one for the spreading. As shown in Figure \ref{fig:fig5}, beam spiraling and axisymmetrization  are observed in both our FSPIC and FSPIF numerical simulations, as expected, in agreement with the results in \cite{cerfon2013analytic, Cerfon_2016}. We use negatively charged particles in all our simulations and hence Figure \ref{fig:fig5} shows the negative of charge density. The same is true for all the rest of the figures where we show charge density from the test cases. In order to demonstrate that the number of modes in Figure \ref{fig:fig5} is sufficient to resolve the spatial scales of this problem, we present in Figure \ref{fig:fig6} the absolute values of the Fourier modes for $N_m=128$, and it is clear that this Fourier resolution is sufficient to capture the part of the Fourier spectrum that has a significant impact on the evolution of the beam.

\begin{figure}
	\centering
   \includegraphics[width=0.8\textwidth]{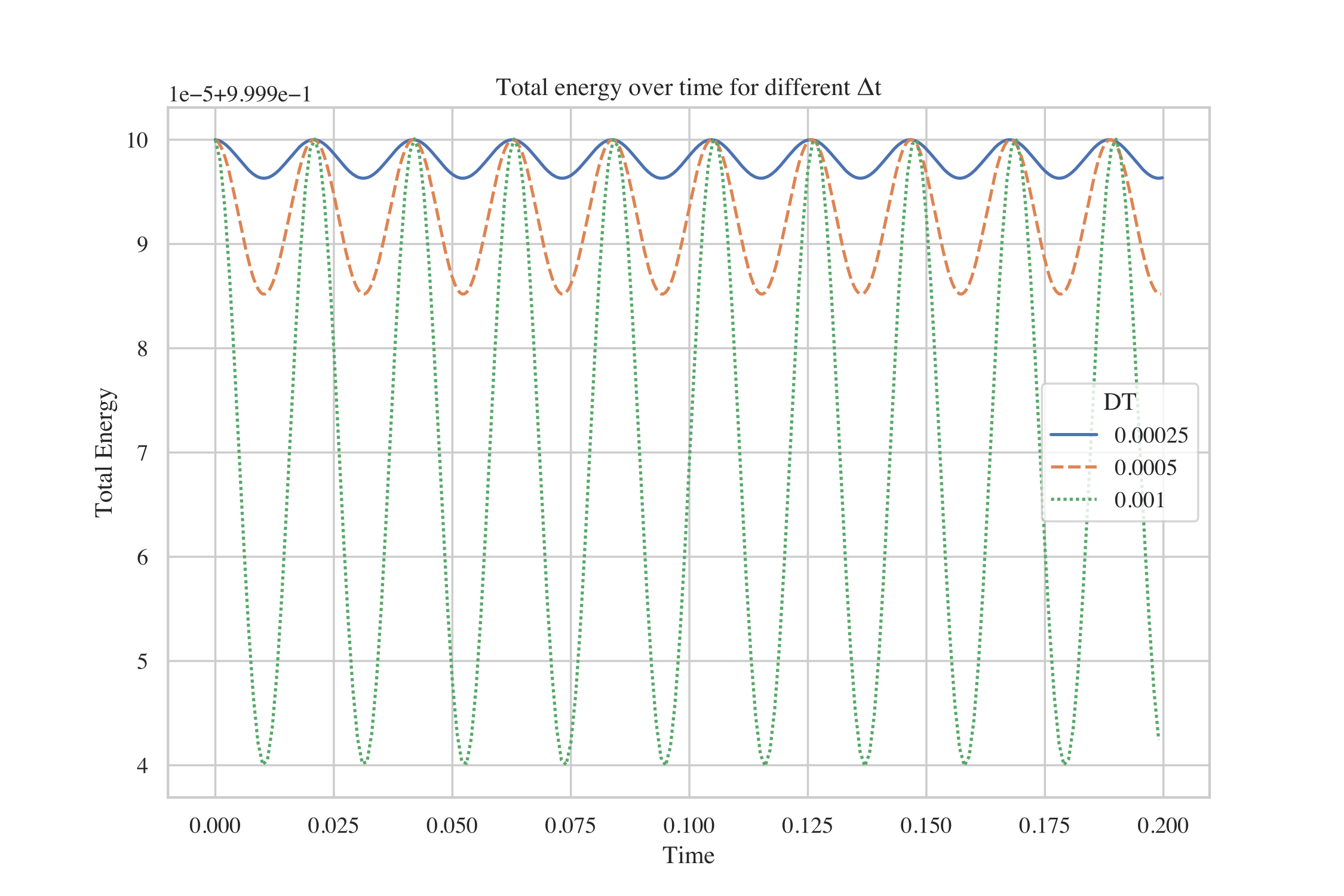}
	\caption{Total energy versus time for the full duration of the simulation and for different time step sizes $\Delta t=0.001, 0.0005,$ and $0.00025$. The periodic variation of the energy has a period equal to the cyclotron period, indicating that the energy fluctuation is mostly caused by the strong magnetic field.}
	\label{fig:fig7}
\end{figure}

\begin{figure}
	\centering
	\includegraphics[width=0.8\textwidth]{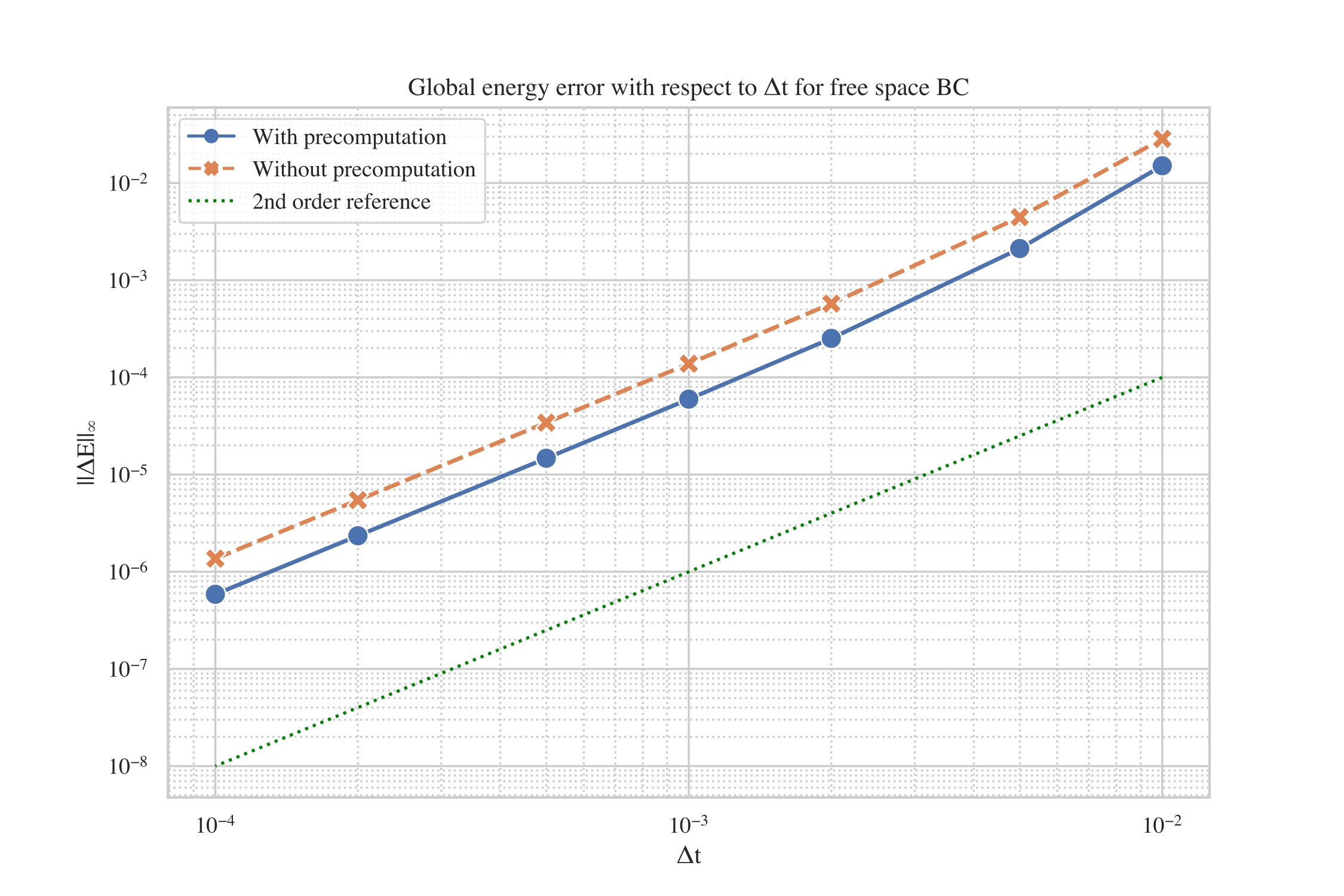}
	\caption{Convergence of the global relative energy error (measured by the $L_\infty$ norm) for different time step sizes $\Delta t$ for our free space particle-in-Fourier simulation using the Boris algorithm. The Fourier resolution is set to $N_m=32$. Second order convergence is observed whether or not we perform a precomputation step. Interestingly, the precomputation step reduces the energy error by a constant compared to the errors without precomputation.}
	\label{fig:fig8}
\end{figure}

We measure the total energy $E$ over the full duration of the simulation for different choices of $\Delta t$. Figure \ref{fig:fig7} shows the energy fluctuation over time for $\Delta t=0.001, 0.0005, 0.00025$. Notice that the fluctuation is smooth and has a period equal to the cyclotron period, which suggests that the error is mainly caused by the strong magnetic field. Since the Boris algorithm is energy-conserving if $\textbf{E}=\textbf{0}$ and at least second order convergent in the presence of an electric field, the global energy error has a second order of convergence, \textit{i.e.}, $\varepsilon(E)\sim \Delta t^2$, as can be seen in Figure \ref{fig:fig8}. Moreover, from Figure \ref{fig:fig8} we observe that even with the precomputation step, the energy still converges quadratically, and in fact the error is smaller than in the case without precomputation by a constant factor.\\

\begin{figure}
	\centering
	\includegraphics[width=0.8\textwidth]{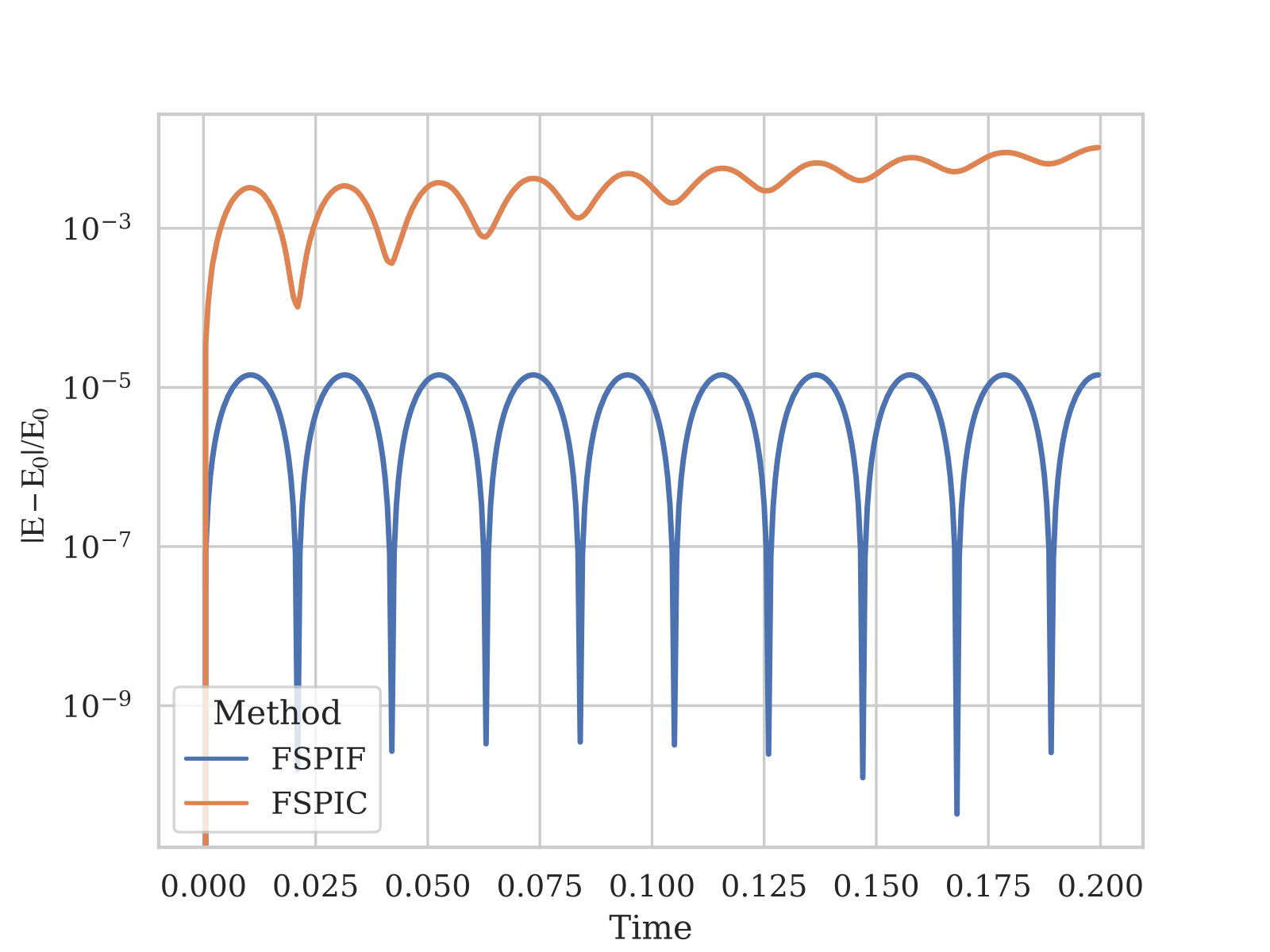}
	\caption{Comparison of the total energy fluctuation as a function of time between our PIC and PIF implementations. Our PIC scheme uses the Vico solver to compute the free space Poisson electrostatic potential on the collocated grid. The time step size is $\Delta t=0.0005$, and the number of modes in PIF and grid points in PIC is set to be $N_g=N_m=32$.}
	\label{fig:fig9}
\end{figure}

We also compare FSPIF and FSPIC with the same number of modes for FSPIF as the number of grid points for FSPIC in each dimension ($N_m=N_g=32$), the same number of particles $N_p=40,000$, and the same time step size $\Delta t=0.0005$, and show the numerical results in Figure \ref{fig:fig9}. Clearly, the use of a grid in FSPIC causes the energy to increase over time, whereas it is very stable with FSPIF. This will be more pronounced for certain long time integration simulations as we show in the next section.

\subsubsection{\myblue{Grid heating in a long time integration simulation}}
\label{grid_heating}

\begin{figure}
\subfigure[FSPIC, $64^2$ grid]{
\includegraphics[width=0.3\columnwidth]{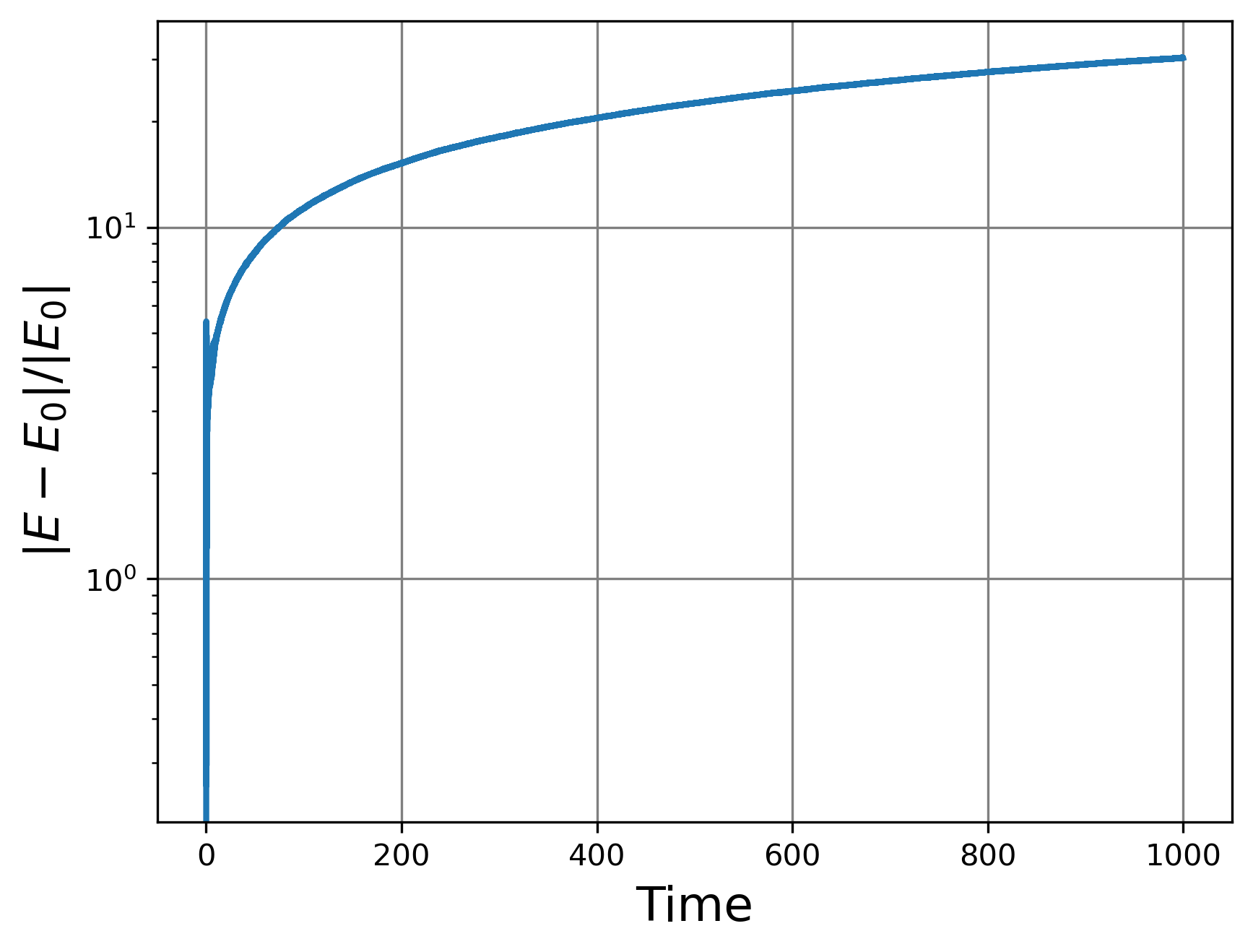}
\label{fig:9apic64}}
\subfigure[FSPIF, $64^2$ modes]{
\includegraphics[width=0.3\columnwidth]{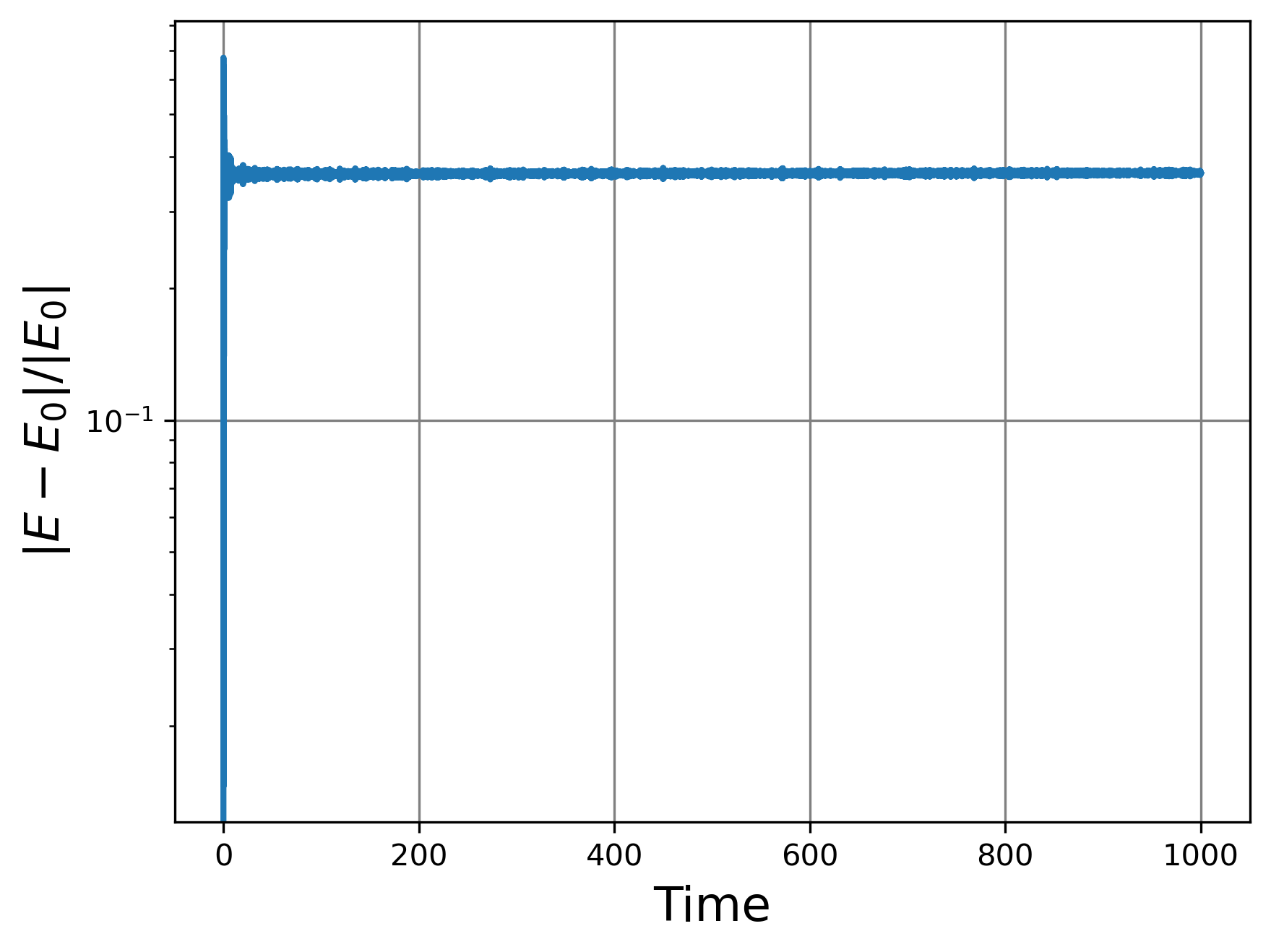}
\label{fig:9apif64}}
\subfigure[FSPIC, $512^2$ grid]{
\includegraphics[width=0.3\columnwidth]{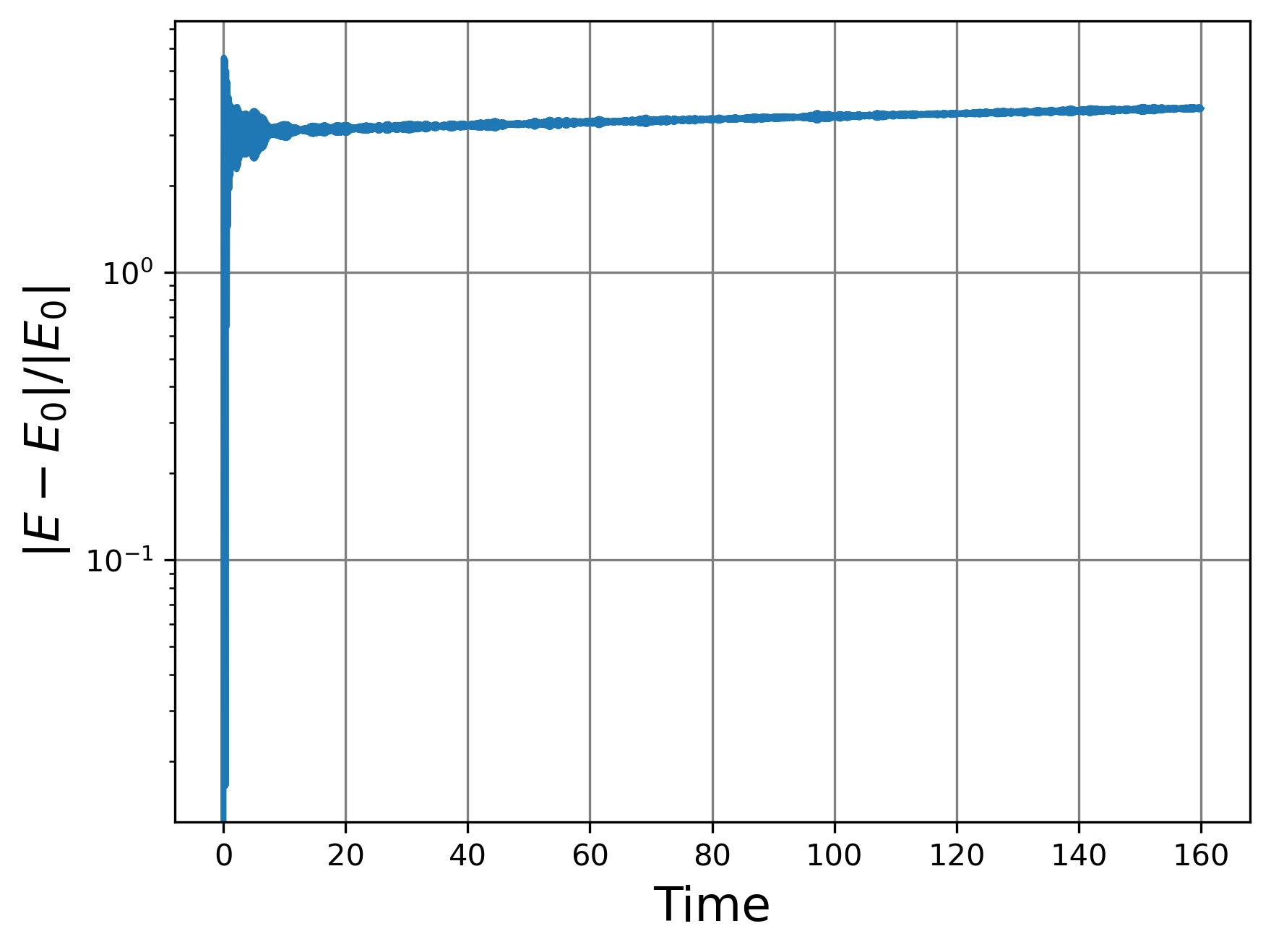}
\label{fig:9apic512}}
 \caption{Comparison of the total energy error as a function of time between FSPIC and FSPIF schemes for the simulations presented in Section \ref{grid_heating}.}
	\label{fig:fig9a}
\end{figure}

\begin{figure}
\subfigure[FSPIC, $64^2$ grid]{
\includegraphics[width=0.3\columnwidth]{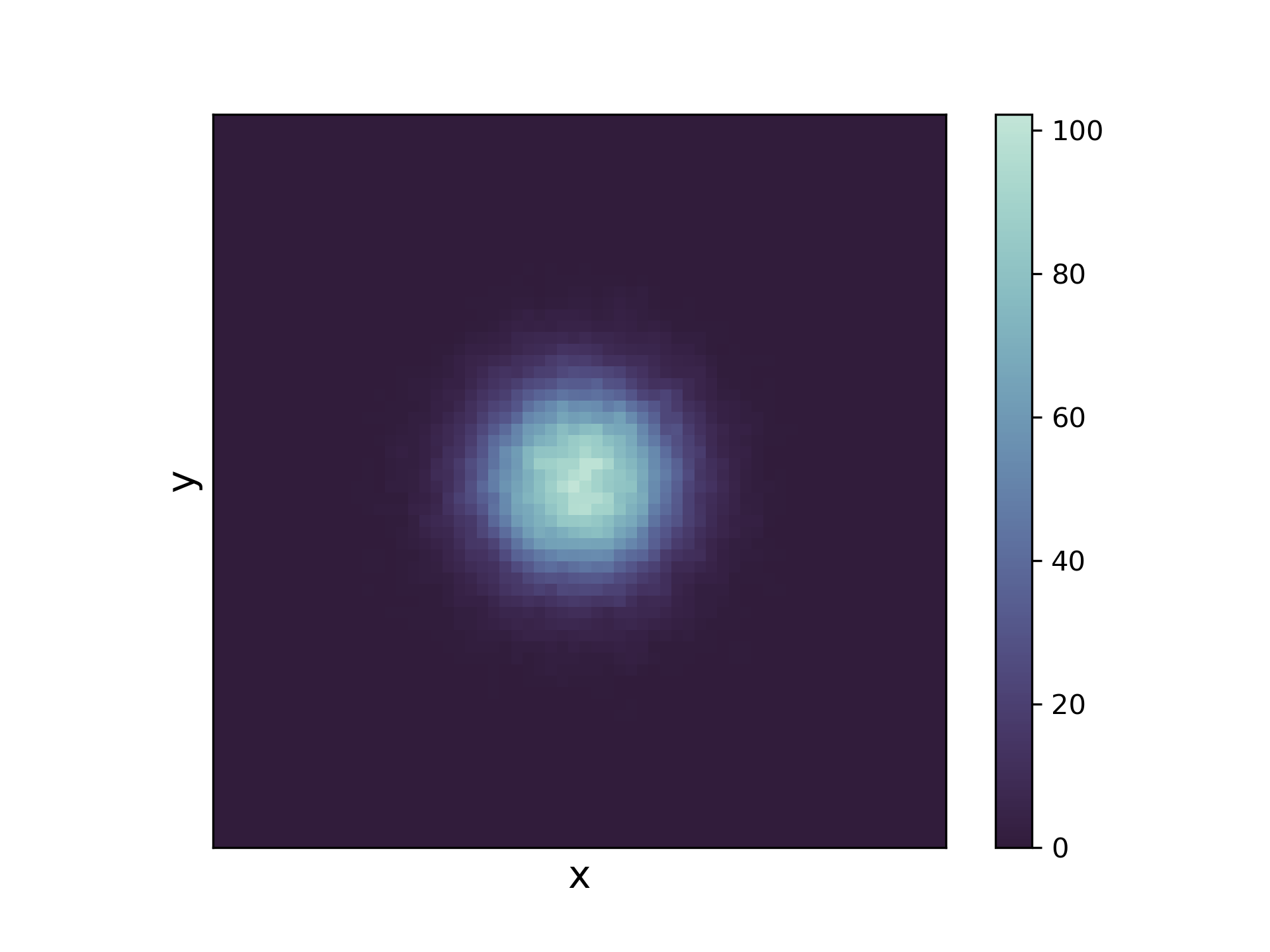}
\label{fig:9bpic64}}
\subfigure[FSPIF, $64^2$ modes]{
\includegraphics[width=0.3\columnwidth]{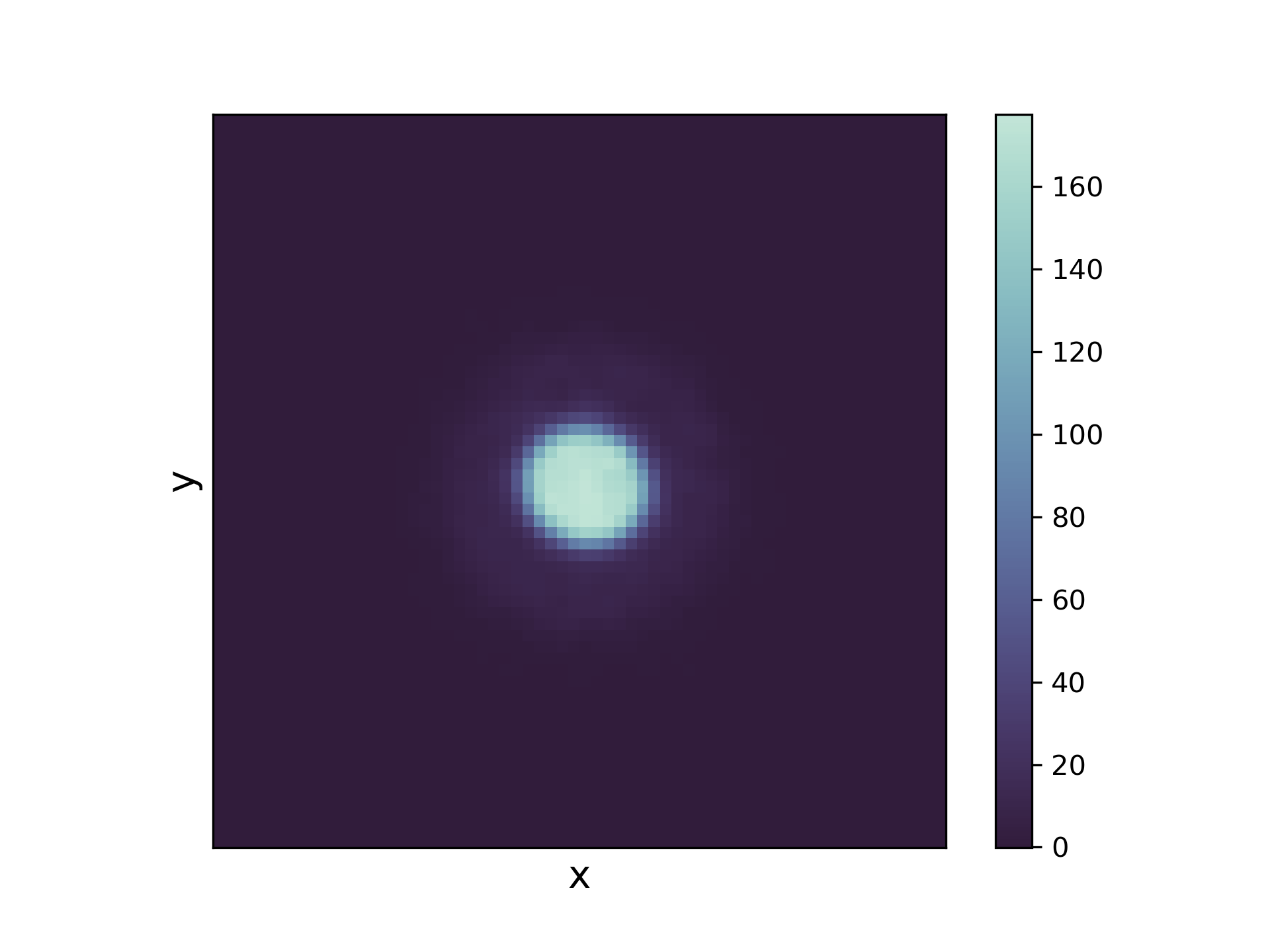}
\label{fig:9bpif64}}
\subfigure[FSPIC, $512^2$ grid]{
\includegraphics[width=0.3\columnwidth]{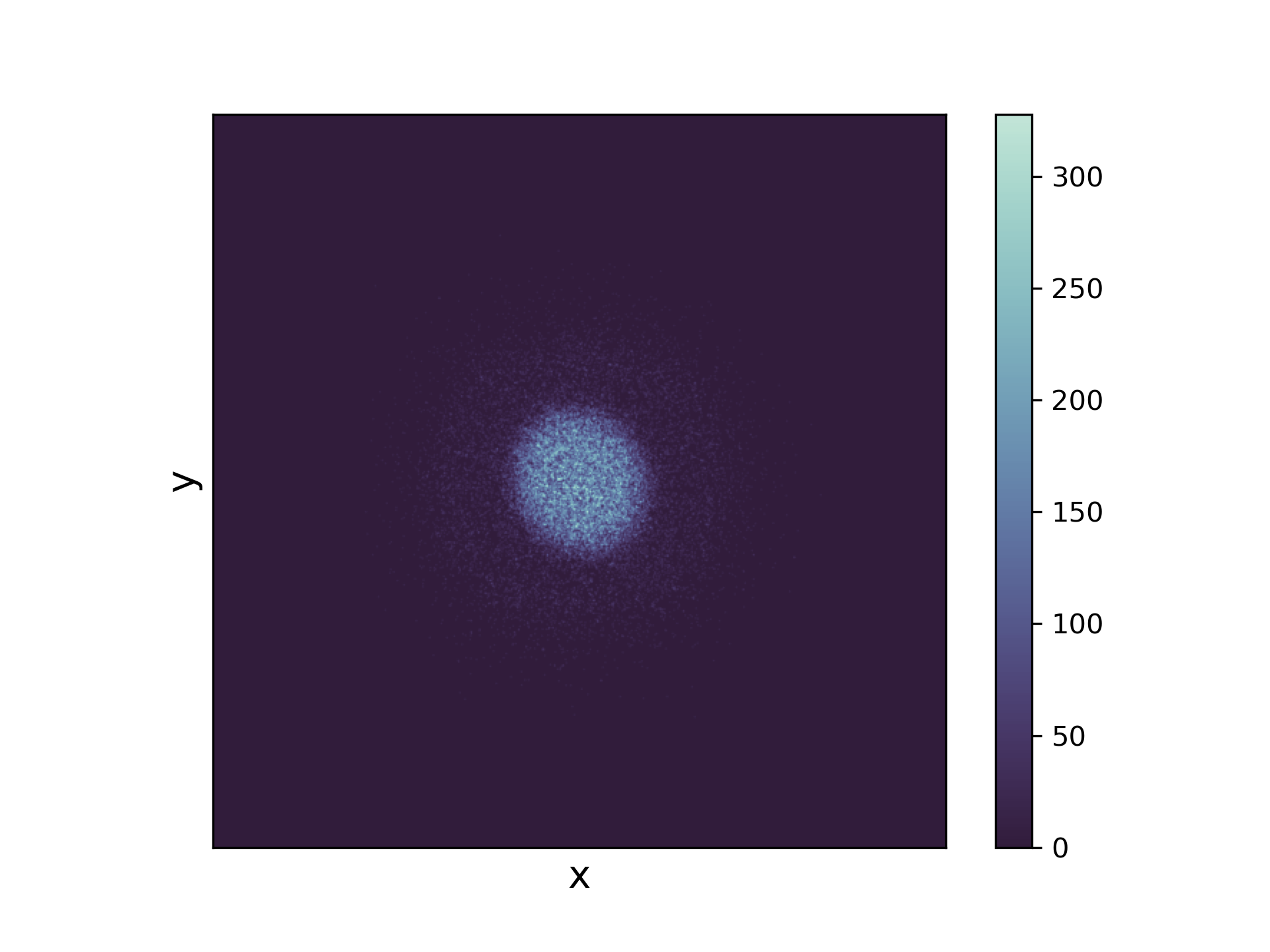}
\label{fig:9bpic512}}
 \caption{Charge density comparison at $T=160$ between FSPIC and FSPIF schemes for the simulations presented in Section \ref{grid_heating}.}
\label{fig:fig9b}
\end{figure}

\begin{figure}
\subfigure[Charge density, FSPIC, $64^2$ grid, $T=1000$]{
\includegraphics[width=0.5\columnwidth]{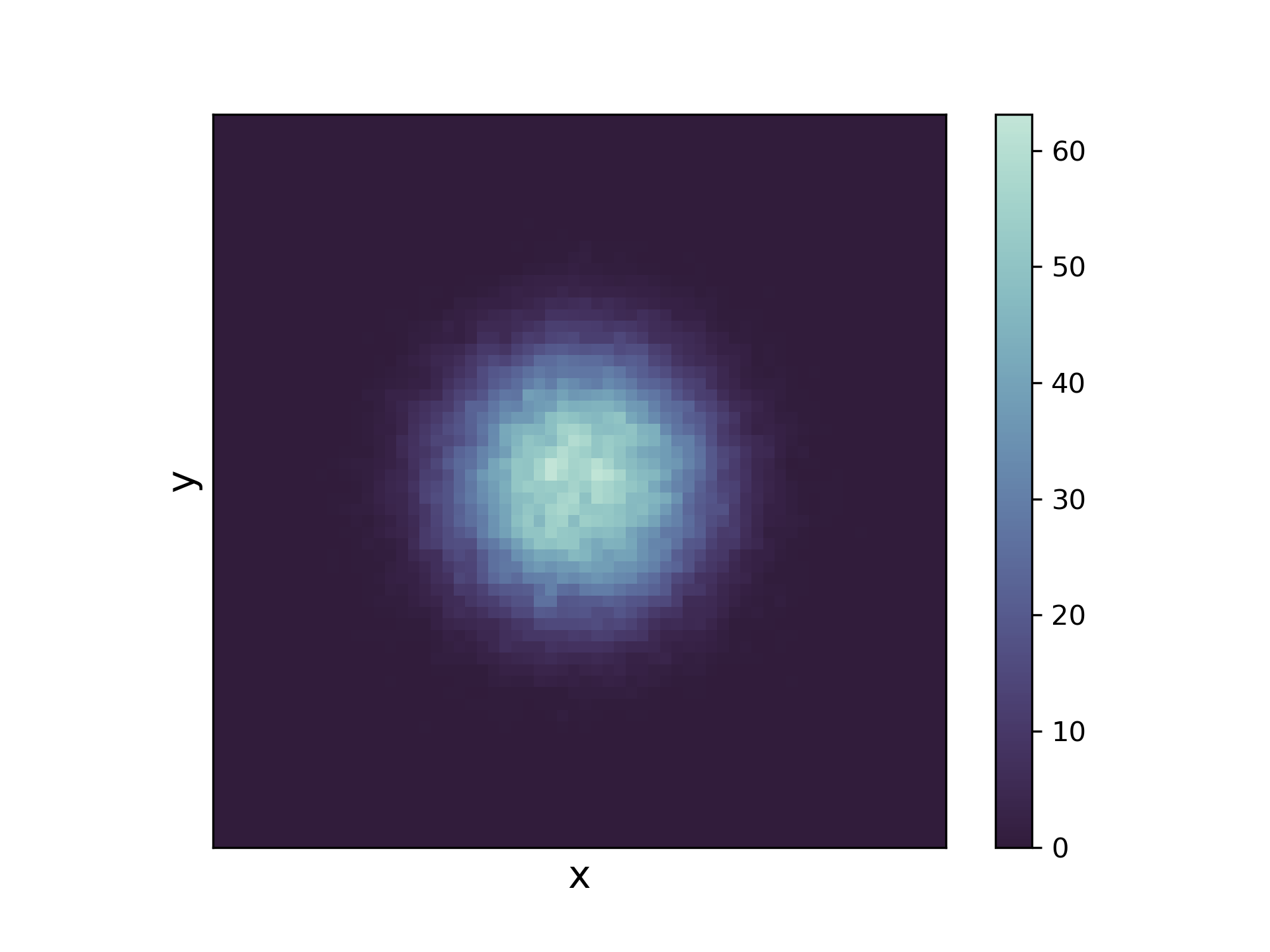}
\label{fig:9cpic64}}
\subfigure[Charge density, FSPIF, $64^2$ modes, $T=1000$]{
\includegraphics[width=0.5\columnwidth]{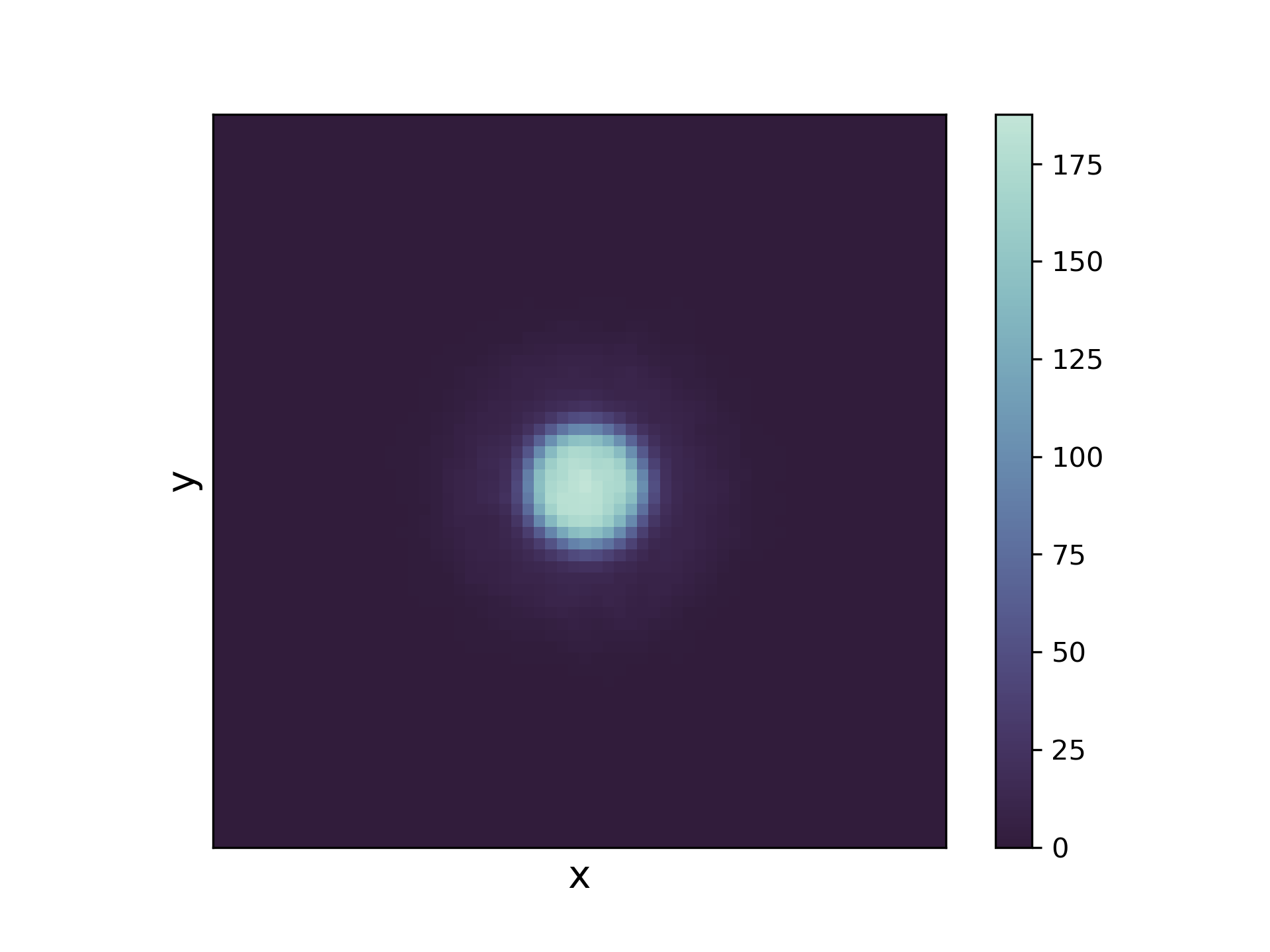}
\label{fig:9cpif64}}
\subfigure[$x-$phase space, $T=160$]{
\includegraphics[width=1\columnwidth]{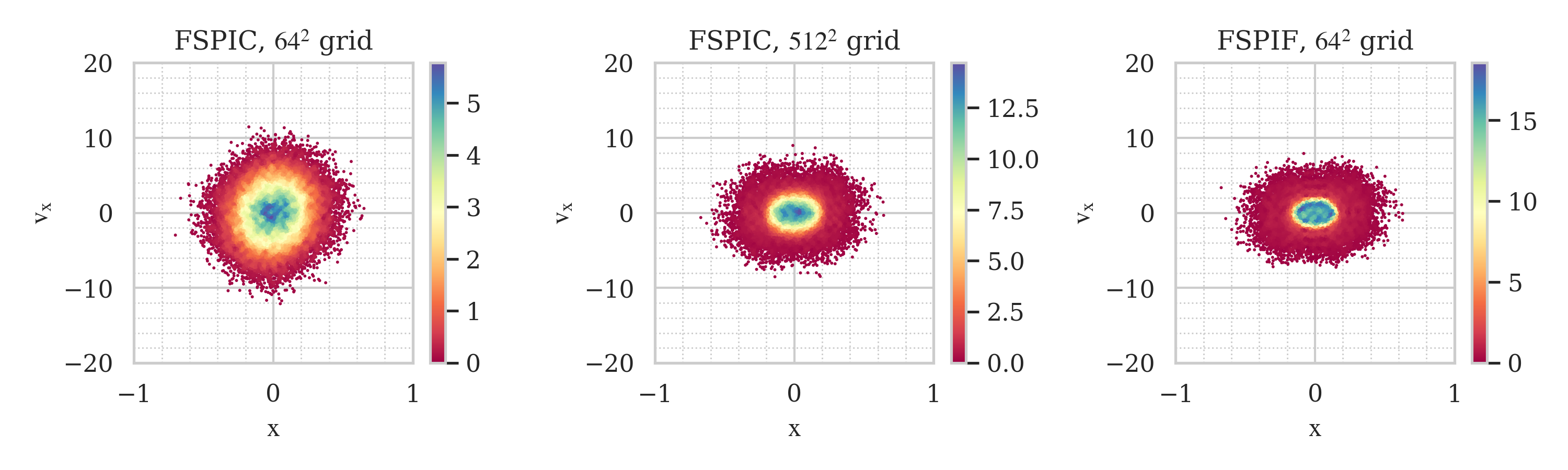}
\label{fig:9cphase160}}
\subfigure[$x-$phase space, $T=1000$]{
\includegraphics[width=1\columnwidth]{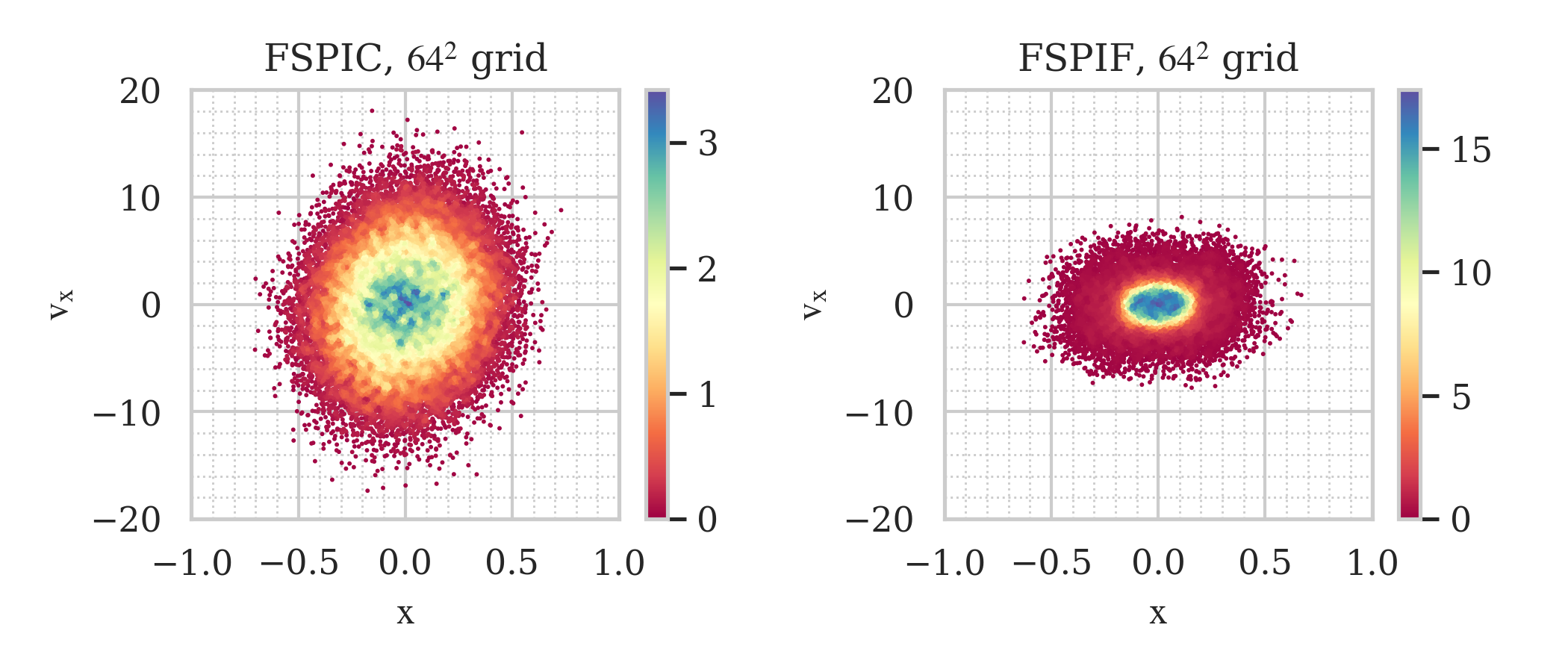}
\label{fig:9cphase1000}}
 \caption{Phase space and charge density comparisons at $T=160$ and $T=1000$ for FSPIC and FSPIF schemes for the simulations presented in Section \ref{grid_heating}.}
\label{fig:fig9c}
\end{figure}

\myblue{In this section, we consider a long time integration simulation, and compare FSPIC and FSPIF schemes in terms of the accuracy of the physical quantities of interest, as well as time to solution. We consider a setting slightly different from the one in Section \ref{cb}, with the domain $\Omega = [-1, 1]\times [-1, 1]$, external magnetic field $\textbf{B}_0=40\mathbf{z}$, total charge $Q=-20$ and the standard deviations in equation \ref{eq:ic_cyclotron} scaled corresponding to the domain as $\sigma_x=1/15$ and $\sigma_y=1/5$. We compare an FSPIF simulation with $N_m=64$ modes with two FSPIC simulations, one with $N_m=64$ grid points in each direction, and one with $N_m=512$ grid points in each direction. We take $\Delta t = 0.01$, total number of particles $40,000$ and integrate until final time $T=1000$ for FSPIC and FSPIF with $N_m=64$, whereas for FSPIC with $N_m=512$ grid points in each direction, we only integrate until $T=160$ to reduce computing time requirements. FSPIF with precomputing is used for the simulations. 

We compare the energy conservation between FSPIC and FSPIF schemes in Figure \ref{fig:fig9a}. While the energy error is very stable for FSPIF with $N_m=64$ in Figure \ref{fig:9apif64}, it is two orders of magnitude larger and increases continuously for FSPIC scheme with $N_m=64$ in Figure \ref{fig:9apic64} due to aliasing and grid heating. For FSPIC with $N_m=512$
the slope of the error increase is much smaller in Figure \ref{fig:9apic512} compared to Figure \ref{fig:9apic64} but the values are still an order of magnitude higher than for the FSPIF scheme in Figure \ref{fig:9apif64}.

The grid heating and lack of energy conservation in FSPIC scheme are reflected physically in the charge density and phase-space plots in Figures \ref{fig:fig9b} and \ref{fig:fig9c}. We can clearly see that the beam radius for the FSPIC scheme with $N_m=64$ is larger in Figure \ref{fig:9bpic64} as compared to the beam radius in Figure \ref{fig:9bpif64} for the FSPIF scheme with $N_m=64$. The beam radius for FSPIC with $N_m=512$ in Figure \ref{fig:9bpic512} is comparable to the FSPIF scheme, however, since the total number of particles is not increased, it has much more statistical noise. We observe similar results in the $x-$phase space plot in Figure \ref{fig:9cphase160} where the $x-$emittance is comparable for FSPIF with $N_m=64$ and FSPIC with $N_m=512$, whereas for FSPIC with $N_m=64$ the semi-minor axis is twice larger. The effects are much more pronounced at $T=1000$ (this corresponds to $10^5$ time steps, which is typical for many production-level cyclotron simulations), where the semi-minor axis is now approximately three times larger in Figure \ref{fig:9cphase1000} and the corresponding beam radius is also larger as shown in Figure \ref{fig:9cpic64}.

 Even though in the earlier tests we used a strict tolerance of $10^{-14}$ close to machine precision  for the NUFFTs in the FSPIF scheme, it can be made lenient depending on the time step size in the explicit time integrator so that the energy conservation is not affected. For the chosen time step size of $\Delta t=0.01$ we found from numerical experiments that the NUFFT tolerance of $10^{-4}$ produces indistinguishable results from that of $10^{-14}$. We note here that while the momentum conservation is unaffected by the tolerance of NUFFT, charge conservation is proportional to it and for this test a charge conservation error of $\mathcal{O}(10^{-4})$ did not change the results in any significant way.

 Now the times to solution for the FSPIC and FSPIF schemes are compared in Table \ref{tab:timepicpif}. The computations are performed in one node of CPU partition of JURECA DC supercomputer at Jülich Supercomputing Centre. Each node of JURECA DC has 128 cores of AMD EPYC 7742, 2.25GHz processors with 512GB RAM. The computations are neither parallelized nor optimized explicitly and the default settings of \texttt{fiNUFFT} as well as other libraries such as numpy and scipy are used. While the FSPIF scheme with $N_m=64$ and NUFFT tolerance of $10^{-4}$ is slightly more expensive than the FSPIC scheme with $N_m=64$ (approximately $1.1$ times), it is approximately $12$ times faster than the FSPIC scheme with $N_m=512$ which is needed for comparable accuracies in physical quantities of interest such as beam radius and emittance. Thus even in 2D-2V, FSPIF gives significant speedup compared to the FSPIC scheme with the fine resolutions necessary to avoid grid heating. We expect much larger benefits in 3D-3V simulations due to the increased number of grid points and particles.} 

 \begin{table}[h!b!t!]
         \centering
         \begin{tabular}{|r|c|c|}
         \hline
             \!\!\! Scheme \!\!\!\! & \!\!\! No. of grid points/modes \!\!\!\! & \!\!\! CPU time (s) \!\!\!\!\\ 
         \hline
                FSPIC & $64^2$ & 272.6 \\
                FSPIF & $64^2$ & 298.6 \\
                FSPIC & $512^2$ & 3557.5 \\
         \hline
         \end{tabular}
         \vspace{5mm}
        \caption{\label{tab:timepicpif} Time to solution for FSPIC and FSPIF schemes with $40,000$ particles, $\Delta t=0.01$ integrated up to $T=160$. FSPIF scheme with $64^2$ modes and FSPIC scheme with $512^2$ grid points have comparable accuracies in quantities of interest however FSPIF is approximately $12$ times faster than FSPIC.}
\end{table}

\subsection{\myred{Expanding plasma in free space}}

\begin{figure}
\subfigure[Charge density, $T=0$]{
\includegraphics[width=0.3\columnwidth]{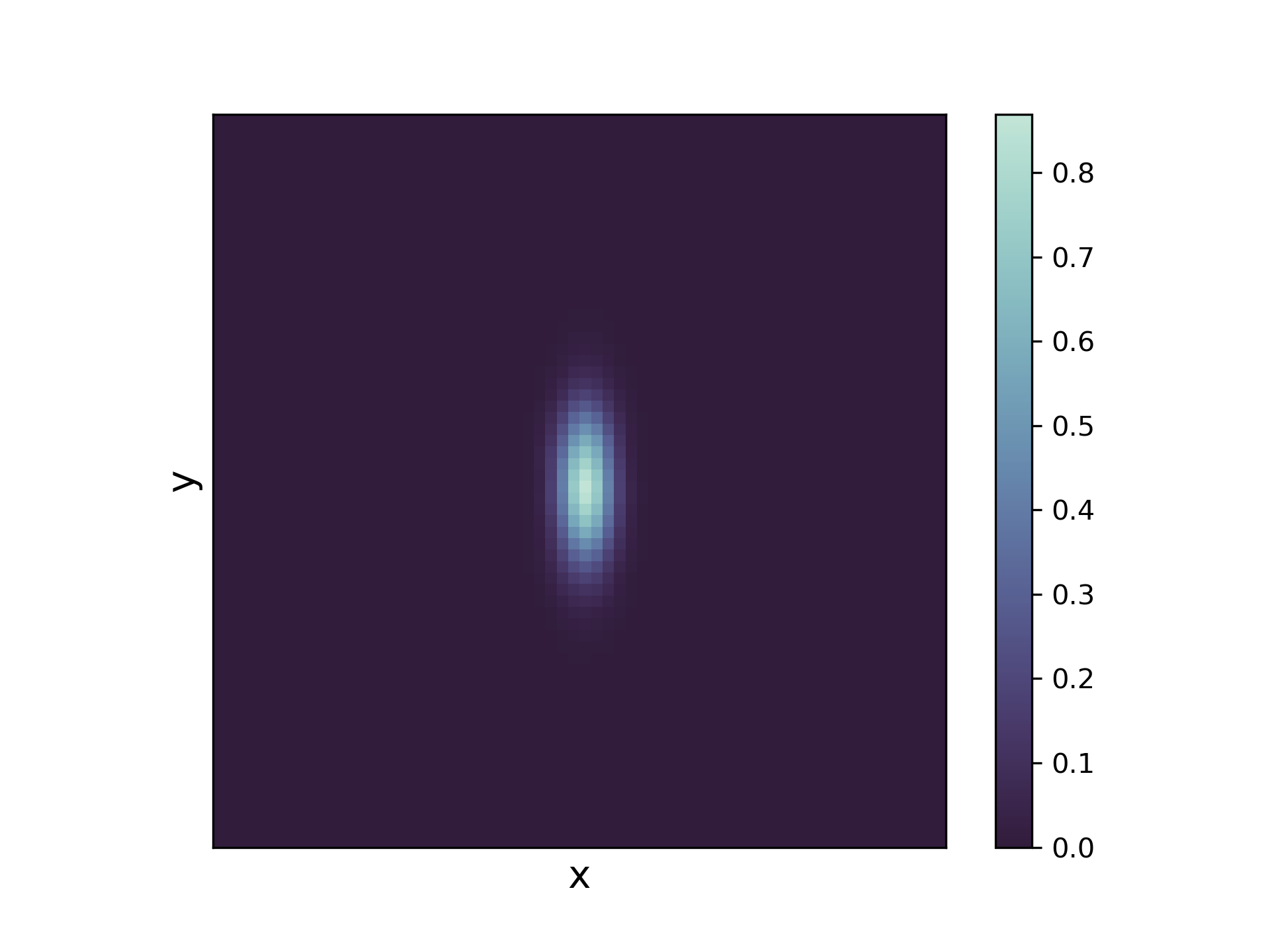}
\label{fig:momconst0}}
\subfigure[Charge density, $T=1.0$]{
\includegraphics[width=0.3\columnwidth]{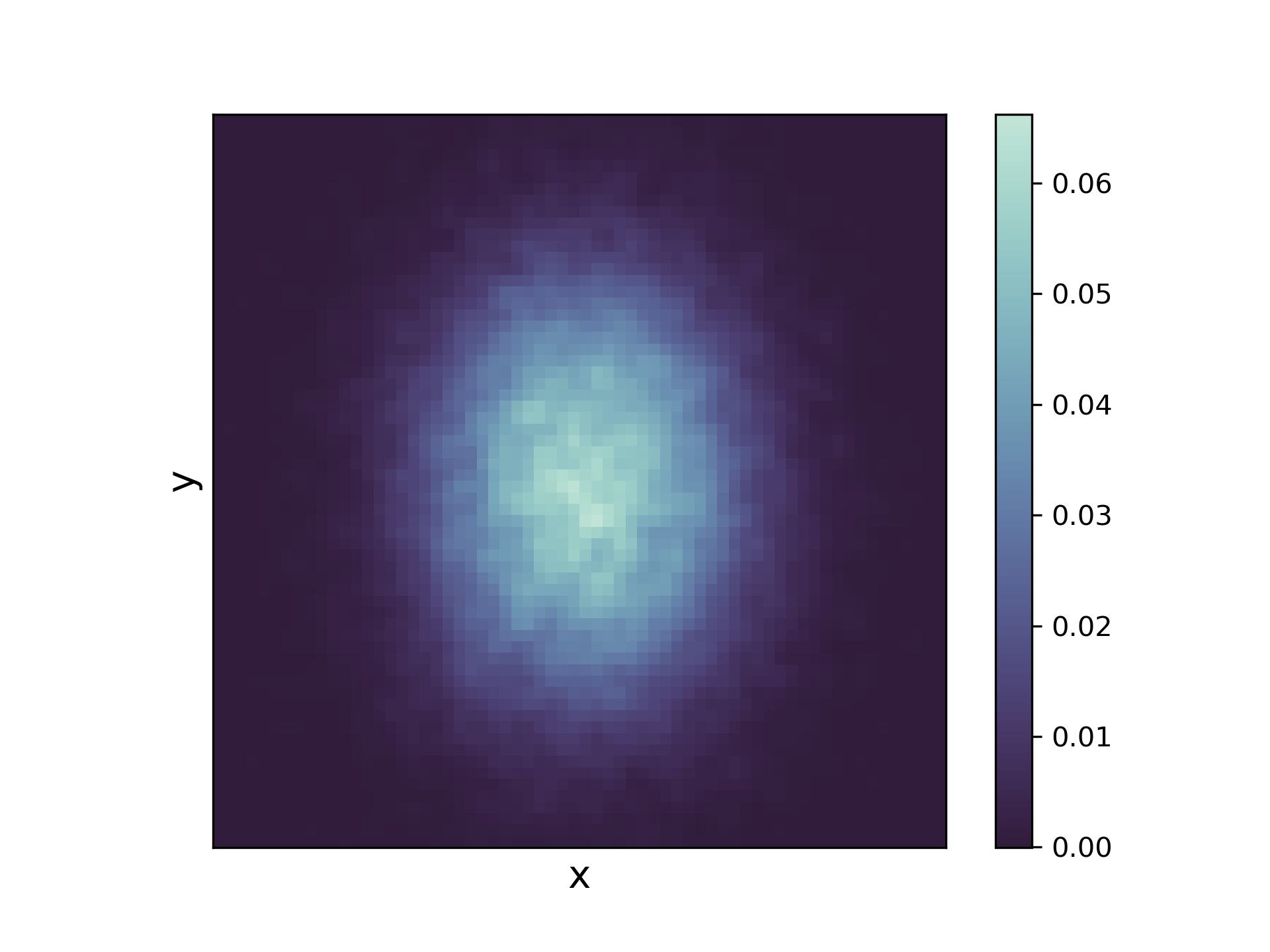}
\label{fig:momconst1}}
\subfigure[Momentum conservation]{
\includegraphics[width=0.3\columnwidth]{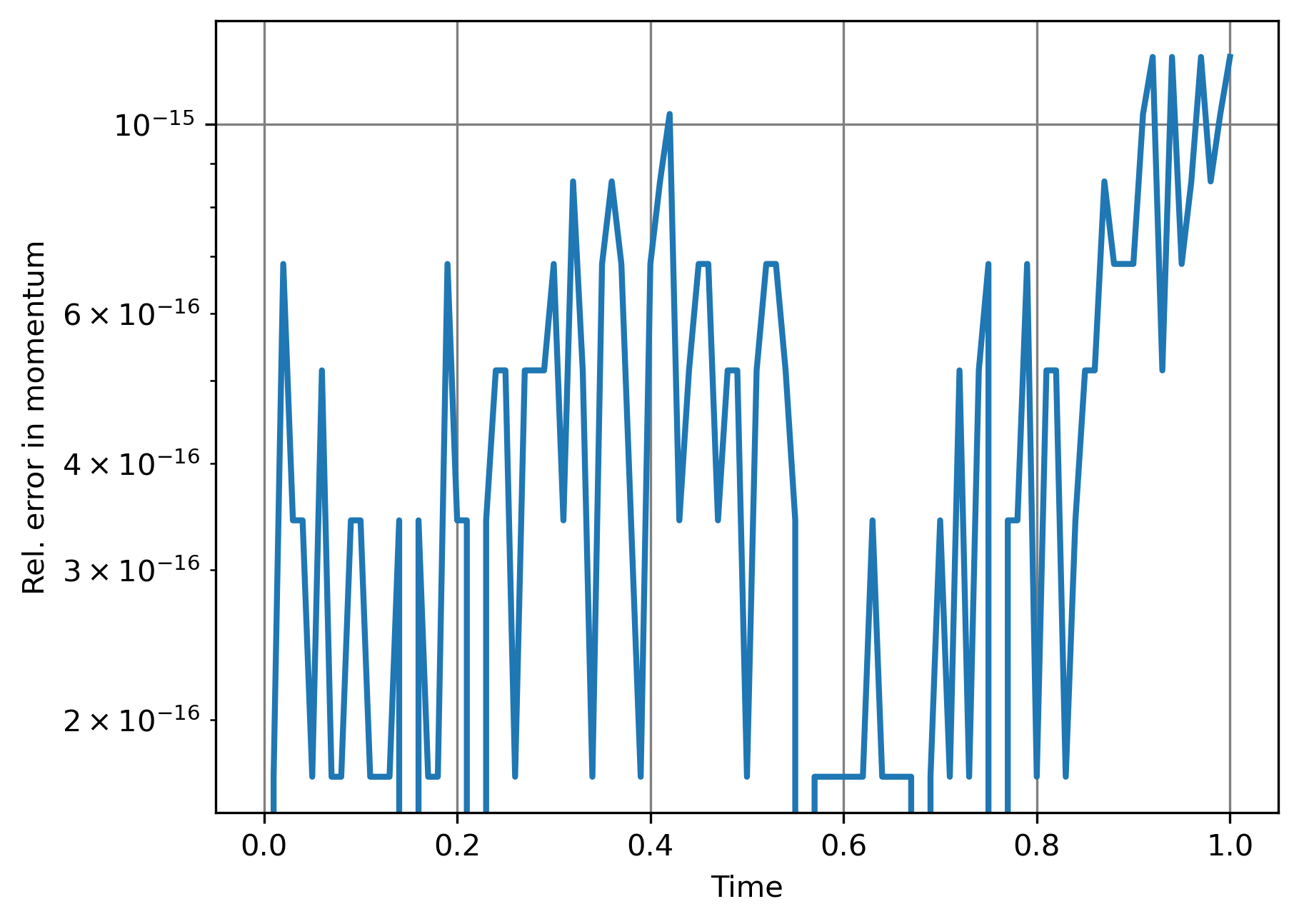}
\label{fig:momconserror}}
 \caption{Charge density at initial time (left), $T=1.0$ (center) and relative error in total momentum conservation (right) for an expanding plasma in free space with the FSPIF scheme.}
\label{fig:momcons}
\end{figure}

\myred{We consider an expanding plasma in free space in this section, and verify the momentum conservation property of the PIF scheme numerically. The parameters for this test case are as follows: domain $\Omega = [-5, 5] \times [-5, 5]$, no external magnetic field, total charge $Q=-1$ and the standard deviations in equation \ref{eq:ic_cyclotron} scaled corresponding to the domain as $\sigma_x=1/3$ and $\sigma_y=1$. We consider a larger domain and a smaller total charge than the previous section so that the expanding plasma does not leave the domain. Similar to the previous section, we take $\Delta t = 0.01$, total number of particles $40,000$, $N_m=64$ and integrate until $T=1.0$. The tolerance for NUFFT is taken as $10^{-4}$ and precomputation is used for the FSPIF scheme. 

In Figure \ref{fig:momcons} we show the initial and final charge densities as well as the relative error in total momentum over time. We get momentum conservation close to double precision accuracy as seen from Figure \ref{fig:momconserror} which confirms our theoretical result in Section \ref{theo_mom_cons}.} 

%
%

\subsection{Infinitely long beam confined by a guiding magnetic field and a perfectly conducting cylinder}
In this section, we consider the same infinitely long beam confined by a guiding magnetic field aligned with it, except that the beam is now also surrounded by a perfectly conducting cylinder. The examples we study serve to verify our implementation of the Dirichlet boundary conditions for our PIF solver, which follows the scheme discussed in Section \ref{abc}. Physically, they are motivated by experiments by Noah Hurst and collaborators \cite{Hurst2016}. 

We begin by first checking the accuracy of our boundary integral solver using the method of manufactured solutions. We choose our analytic solution to be a harmonic function \begin{equation}
    u(x,y) = x^2-y^2 + 2.
\end{equation} Let the coordinate of a boundary point be described by a single variable $\theta$ where $\sin(\theta)=x$ and $\cos(\theta)=y$, then the boundary data for our Laplace equation is:
\begin{equation}
    f(\theta) = \sin^2(\theta)-\cos^2(\theta) + 2.
\end{equation}
\begin{figure}
	\centering
	\includegraphics[width=1\textwidth]{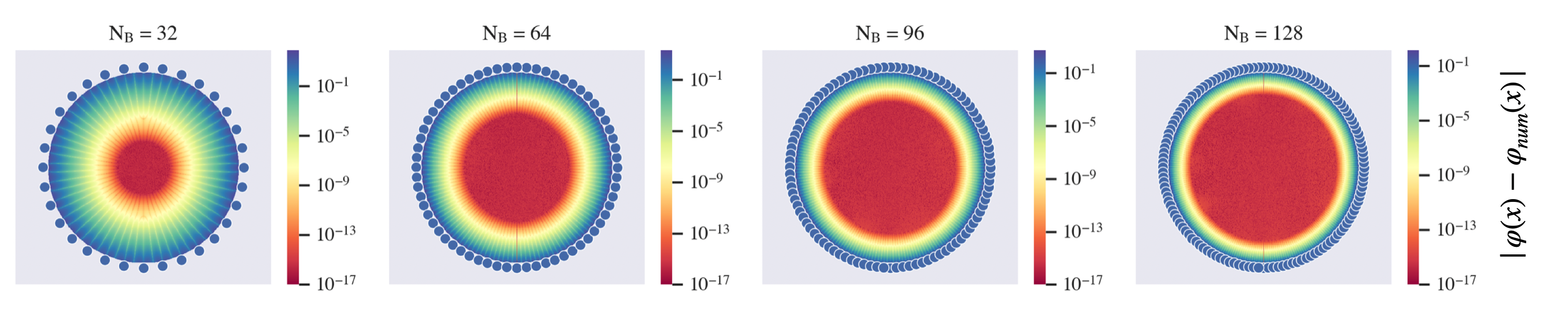}
	\caption{Convergence of the boundary integral method as a function of the number of boundary points $N_B$ in the domain $D_{1}(\V{0})$. The method is spectrally accurate inside the domain except for the regions in the vicinity of the boundary.}
	\label{fig:fig10}
\end{figure}
We vary the number of boundary points and use the trapezoidal rule to evaluate equation \eqref{eq:nlap} on a $256\times 256$ equi-spaced grid that covers the entire domain $D_1(\V{0})$. Due to the singular integrand, the numerical value $u_{\text{num}}$ does not converge in the vicinity of the boundary. This behavior is expected, and advanced quadrature schemes have been developed to address this issue \cite{klockner2013quadrature}. However, the implementation of such schemes is not required if the beam remains sufficiently far from the conductor during its evolution, and is beyond the scope of this work. We therefore only measure the error $\varepsilon(\V{x})=|u_{\text{num}}-u|(\V{x})$ within a slightly shrunk domain $D_{0.9}(\V{0})$, and plot them in Figure \ref{fig:fig10}. Overall, the solver provides spectral accuracy inside $D_{0.9}(\V{0})$ away from the boundary, which matches the spectral convergence of the 
trapezoidal rule for periodic analytic functions \cite{trefethen2000spectral, trefethen2014exponentially}.

\begin{figure}
	\centering
	\includegraphics[width=1\textwidth]{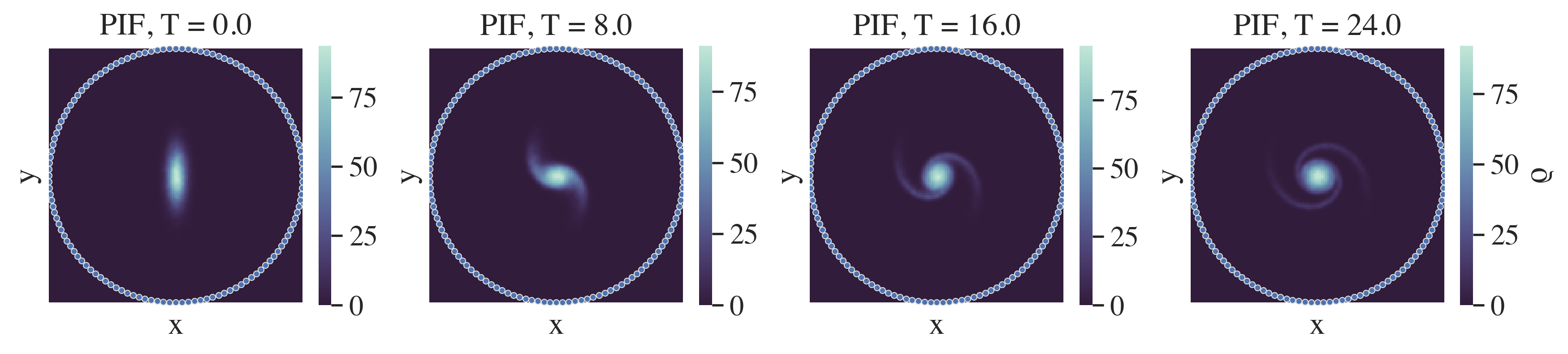}
	\caption{Charge density $\rho$ at different times of an infinitely long non-neutral beam confined by a strong constant and uniform magnetic field and a grounded cylinder. We choose 128 boundary points to evaluate $\varphi^H$ and $\V{E}^H$ using the method described in Section \ref{abc}. The number of modes is $N_m=128$, and the number of particles is 40,000. The solution is similar to the vacuum solution in Figure \ref{fig:fig5}, as expected.}
	\label{fig:fig11}
\end{figure}

\begin{figure}
	\centering
	\includegraphics[width=1\textwidth]{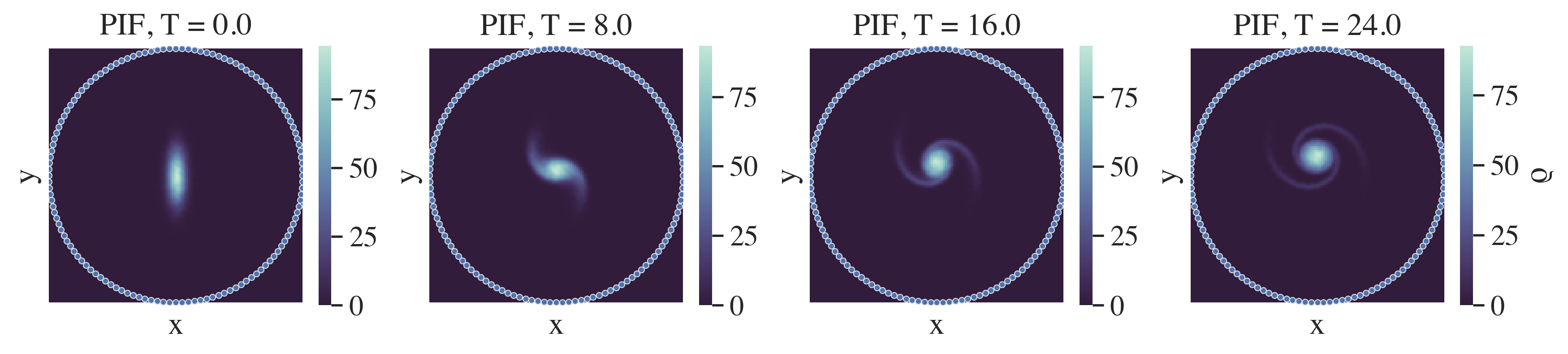}
	\caption{Charge density $\rho$ at different times of an infinitely long non-neutral beam confined by a strong constant and uniform magnetic field and a cylinder with biasing potential $f_2(\textbf{x}) = y$ . We choose 128 boundary points to evaluate $\varphi^H$ and $\V{E}^H$ using the method described in Section \ref{abc}. We choose the number of modes to be $N_m=128$ and the number of particles to be $40,000$. The beam drifts upwards in the $y$ direction over time, as expected.}
	\label{fig:fig12}
\end{figure}

Now we combine our boundary integral solver with the free space PIF algorithm as described in Section \ref{abc}, and test its numerical properties. We consider the same initial charge distribution as in
Section \ref{cb}. In order to keep our computation in the unit box, in this test we use $D_{0.5}(\textbf{0})$ as our domain, instead of the unit disk as in Section \ref{abc}. When solving for the fields and accelerations, we add the harmonic function from the boundary integral method in order to satisfy Dirichlet boundary conditions at the boundary, $C_{0.5}(\V{0})$. We perform tests for two Dirichlet problems, described as 
\begin{align}
    f_1(\textbf{x}) = 0 \quad \forall x \in C_{0.5} (\textbf{0}), \\ f_2(\textbf{x}) = y \quad \forall x \in C_{0.5} (\textbf{0}).
\end{align}
The evolution of the system under these two different Dirichlet boundary conditions over time is plotted in Figures \ref{fig:fig11} and \ref{fig:fig12} with Fourier resolution $N_m=128$, where the blue scatter dots are the 128 boundary points where we evaluate the boundary condition. For the first Dirichlet BC problem $f=f_1$, we expect the solution to be similar to the free-space solution; for the second Dirichlet BC problem $f=f_2$, we expect the beam to receive an extra force along the $y$ axis, so the non-neutral beam will drift upward in the $y$ direction.

\begin{figure}
	\centering
     \includegraphics[width=0.8\textwidth]{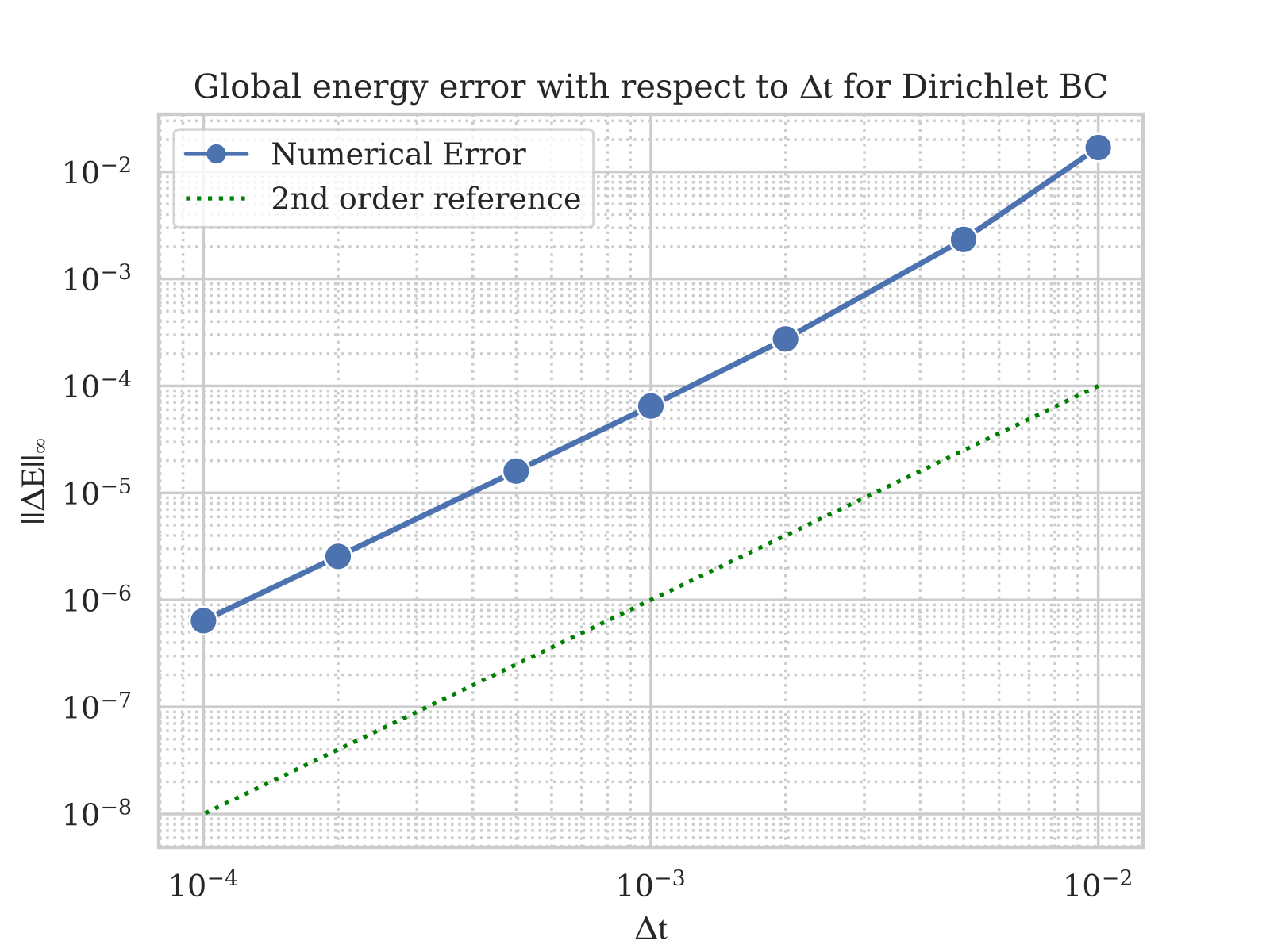}
	\caption{Convergence of the global relative energy error for different time step sizes $\Delta t$ for our particle-in-Fourier simulation with Dirichlet boundary conditions $f_1(\textbf{x})$, corresponding to a grounded, perfectly conducting cylinder. Here, the time integrator is the Boris algorithm. The Fourier resolution is set to $N_m=32$. As for the free space version of our scheme, second order convergence is observed.}
	\label{fig:fig13}
\end{figure}

Our simulations, shown in Figures \ref{fig:fig11} and \ref{fig:fig12}, confirm our intuition. The global relative energy error measured with the $L_\infty$ norm is plotted as a function of the time step size $\Delta t$ in Figure \ref{fig:fig13}. As in the free space case, the energy error converges in second order with respect to $\Delta t$.
\section{Conclusion}\label{sec:con}
We present a novel Particle-In-Fourier (PIF) scheme which can be applied to simulations with free space boundary conditions while maintaining semi-discrete energy conservation. We show that, by using a mollified Green's function kernel, we can obtain highly accurate solutions to the free space Poisson equation, even for non-equispaced data. Combining this Poisson solver with the PIF scheme allows us to simulate the evolution of plasmas described by the Vlasov equation with free space boundary conditions. We prove the energy and momentum conservation of the scheme in semi-discrete and fully discrete settings and validate it with numerical results. 

Our approach can be extended to general Dirichlet boundary conditions on general domains, via standard potential theory. As an illustration of this, we presented an algorithm combining the free space PIF scheme with the boundary integral method on a unit disk for general Dirichlet boundary conditions. We show similar accuracy and energy conservation properties as those of our free space PIF scheme. Numerical implementation and validation of three-dimensional examples in arbitrary geometries are left as future work.

The high accuracy of the field solver and conservation of energy observed in our numerical tests show potential and advantages for long-time three-dimensional plasma simulations based on particle codes. It is however important to remember that in most cases, the error in a particle code is dominated by the noise of particle sampling rather than numerical error due to the field solver. Therefore, in order to further increase the accuracy and performance of our scheme, noise reduction strategies are of particular interest. Noise reduction is typically achieved by running simulations with a larger number of computational particles, taking advantage of scalable massively parallel computing. The design of efficient parallelization strategies for our PIF schemes therefore represents a natural direction for future work. In a complementary way, sparse grids have been used to enhance the performance of PIC by drastically reducing particle sampling error for a given number of computational particles \cite{ricketson2016sparse, muralikrishnan2021sparse}. The development of counterparts to sparse grid PIC schemes in the PIF setting appears to be another promising research direction. 

\section*{Acknowledgements}
C.N.S. was partially supported by National Science Foundation (award DMS-1646339) via the Applied Math Summer Undergraduate Research Experience at NYU Courant. A.J.C. was partially supported by the United States National Science Foundation under Grant
No. PHY-1820852. We would like to acknowledge the computing time in the JURECA DC supercomputer provided under the project CSTMA.
\section*{Appendix A. Radially Symmetric Shape Functions}\label{rssf}
In this work we use radially symmetric shape functions in our PIF code to conduct numerical experiments. These shape functions preserve exact rotational invariance as compared to the tensor-product shape functions used in PIC. The zeroth order radially symmetric b-spline function in 2D can be defined as \begin{equation}
    S_0(\textbf{x}) = \begin{cases}
        \frac{4}{\pi} & \text{if}\ |\textbf{x}| < 1/2,\ \textbf{x}\in\mathbb{R}^2\\
        \\
        0 & \text{otherwise}.
    \end{cases}
\end{equation}
The $l$-th order radially symmetric b-spline function is defined as \begin{equation}
    S_l(\textbf{x}) = \int_{\mathbb{R}^2} S_0(\textbf{x}-\textbf{x}')S_{l-1}(\textbf{x}') d\textbf{x}'.
\end{equation}
The Fourier modes of the $l$-th order radially symmetric b-spline function in 2D is calculated as:
\begin{equation}
    \hat{S}_l (\textbf{k}) = \left(\frac{J_1 (k)}{k}\right)^l
\end{equation} where $k=|\textbf{k}|$.

In addition, in the free space Poisson solver test (Section \ref{fspp}), we use a spatially truncated Gaussian function as our shape function:\begin{equation}
    S^{R}(\textbf{x}; \sigma) = \begin{cases}
        \frac{1}{\sqrt{2\pi \sigma^2}} \exp\left(-\frac{||\textbf{x}||^2}{2\sigma^2}\right) & \text{if}\ ||\textbf{x}||<R\\ \\
        0 & \text{otherwise}.
    \end{cases}
\end{equation}
The Fourier transform of the above shape function is:\begin{equation}
    \hat{S}^{R}(\textbf{k}; \sigma) = \exp\left(-\frac{\sigma^2 k^2}{2}\right)\Re\left(\text{erf}\left(\frac{R}{\sqrt{2}\sigma}+\frac{i\textbf{k}\sigma}{\sqrt{2}}\right)\right).
\end{equation}

\section*{References}

\bibliographystyle{elsarticle-num}

\bibliography{references}

\end{document}